\documentclass[12pt]{article}

\usepackage{amssymb}

\def\M{{\cal M}}
\def\MM{{\overline{\cal M}}}
\def\ocC{{\overline{\cal C}}}
\def\ocV{{\overline{\cal V}}}
\def\ocW{{\overline{\cal W}}}
\def\R{{\mathbb R}}
\def\C{{\mathbb C}}
\def\L{{\cal L}}
\def\tphi{{\widetilde \varphi}}
\def\CP{{\mathbb C} {\rm P}}
\def\Z{{\mathbb Z}}
\def\Aut{{\rm Aut}}

\def\T{{\widehat{\cal T}}}
\def\B{{\cal B}}
\def\qed{{\hfill $\diamond$}}
\def\d{{\partial}}
\def\Stab{{\rm Stab}}
\def\Gr{{\rm Gr}}

\usepackage{theorem}

\newtheorem{theorem}{Theorem}
\newtheorem{proposition}{Proposition}[section]

{\theorembodyfont{\rmfamily}

\newtheorem{definition}[proposition]{Definition}
\newtheorem{example}[proposition]{Example}
\newtheorem{remark}[proposition]{Remark}
}

\include{epsf}

\title{Strebel differentials on stable curves and Kontsevich's
proof of Witten's conjecture}

\author{Dimitri Zvonkine\thanks{
Institut f\"ur Mathematik,
Universit\"at Z\"urich,
Winterthurerstrasse 190,
CH-8057 Z\"urich, Switzerland.
e-mail: zvonkine@math.unizh.ch}}

\date{\today}

\begin{document}

\maketitle

\begin{abstract}

We define Strebel differentials for stable complex curves,
prove the existence and uniqueness theorem that generalizes
Strebel's theorem for smooth curves, prove that Strebel 
differentials form a continuous family over the
moduli space of stable curves, and show how this 
construction can be applied to clarify a delicate point
in Kontsevich's proof of Witten's conjecture.

\end{abstract}

\section{Introduction} \label{Sec:intro}

\subsection{Motivation}

The main motivation of this paper is to clarify a
delicate point in Kontsevich's proof of the Witten 
conjecture~\cite{Kontsevich}. The conjecture concerns
the intersection numbers of the first Chern classes
of some line bundles $\L_1, \dots, \L_n$ over
the Deligne-Mumford compactification $\MM_{g,n}$
of the moduli space of $n$-pointed genus $g$ curves.

The proof by Kontsevich uses in an essential way the cell 
decomposition of the space $\M_{g,n} \times \R_+^n$
given by Strebel differentials. This cell decomposition
has a natural closure, and it is necessary to know the 
relation between this closure and
the Deligne-Mumford compactification of the moduli
space. In his paper Kontsevich gives the answer, but
the proof is only briefly sketched. 

A full proof was provided in a very 
thorough paper by E.~Looijenga~\cite{Looijenga95}
(and here we give a new proof).
However many people are not sure what the
exact implications of Looijenga's results are.
In particular, S.~P.~Novikov
still insists that there is a gap in
Kontsevich's proof. Therefore we think that it is useful
to give a complete account of the situation. This
paper is essentially an overview, although it
contains some new results 
(Theorems~\ref{Thm:contfam} and~\ref{Thm:cont2} as well as
a new proof of Theorem~\ref{Thm:homeomorphism}).

We assume that the notion of moduli space $\M_{g,n}$
of Riemann surfaces of genus $g$ with $n$ marked
and numbered points is known, as well as the notion
of a stable curve and that of the Deligne-Mumford compactification
$\MM_{g,n}$ of the moduli space. We also use
the standard notation $\L_i$ for the line bundle over
$\MM_{g,n}$ whose fiber over a point $x \in \MM_{g,n}$
representing a stable curve $C_x$ is the cotangent line to $C_x$
at the $i$th marked point.

\subsection{The main steps of the argument}

Here is a sketch of what is to be or has been done to
justify Kontsevich's expressions of the first
Chern classes of the bundles $\L_i$.

{\bf 1.} A quotient $K\MM_{g,n}$
of the Deligne-Mumford compactification $\MM_{g,n}$
by some equivalence relation is constructed. The
quotient is a compact Hausdorff topological orbifold.
The line bundles $\L_i$ on $\MM_{g,n}$ are pull-backs
of some line bundles over $K\MM_{g,n}$, which we will
also denote by $\L_i$. The quotient $K\MM_{g,n}$
has only singularities of real codimension at least~2. Therefore
its fundamental homology class is well-defined.
The fundamental homology class of $\MM_{g,n}$
is sent to the fundamental homology class of $K\MM_{g,n}$
under the factorization.

All these facts are quite simple once the equivalence
relation is given. They are formulated in Kontsevich's
original paper, except for the fact
that $K\MM_{g,n}$ is Hausdorff, which was proved
by Looijenga.

{\bf 2.} $K\MM_{g,n} \times \R_+^n$ is homeomorphic to
a cell complex $A$. There is a piecewise affine projection
from $A$ to $\R_+^n$ that commutes with the projection
from $K\MM_{g,n} \times \R_+^n$ to $\R_+^n$.

The cell complex and the homeomorphism
were constructed by Kontsevich, but he omitted
the proof of continuity. The proof is one of
the main results of Looijenga's paper. Here we
give a different proof.

{\bf 3.} For each line bundle $\L_i$, one
constructs a cell complex $\B_i$ together with a projection
$\B_i \rightarrow A$ satisfying the following properties.
Each fiber of the projection is homeomorphic to a circle.
The (real) spherization 
$(\L_i^* \setminus \mbox{zero section})/\R_+$ of the dual line bundle 
$\L_i^*$ is isomorphic to the circle
bundle $\B_i$.

The bundles $\B_i$ were constructed by Kontsevich, the isomorphism
of the two bundles is immediate.

{\bf 4.} A connection on the bundle $\B_i$ can be explicitly
given. Its curvature represents the first Chern
class of the bundle.

The connection and the curvature were written out
by Kontsevich. However, they are only piecewise
smooth forms defined on a cell complex that is
not homeomorphic to a smooth orbifold. Therefore
one should explain why (and in what sense)
the curvature indeed
represents the first Chern class correctly. This
issue does not seem to have been raised in the
literature. We deal with it here in Section~\ref{Sec:Chern}.

Thus the goal of this paper is to give a complete proof of
the following theorem.

Consider a point 
${\bf p} = (p_1, \dots, p_n ) \in \R_+^n$. Denote
by $A_{\bf p}$ the preimage of $\bf p$ under the
projection from $A$ to $\R_+^n$. Then $A_{\bf p}$
is also a cell complex.

\begin{theorem} \label{Thm:int=int}
To each bundle $\L_i$ one can assign a piecewise
smooth $2$-form $\omega_i$ on the cell complex $A$
in such a way that
$$
\int\limits_{\MM_{g,n}} c_1(\L_1)^{d_1} \dots c_1(\L_n)^{d_n}
\;
=
\int\limits_{
\begin{array}{c}
\mbox{\rm cells of top}\\
\mbox{\rm dimension of } A_{\bf p}
\end{array}
}
\omega_1^{d_1} \dots \omega_n^{d_n}
$$
for every ${\bf p} \in \R_+^n$.
\end{theorem}

\subsection{The organization of the paper}

This paper is organized as follows.

In Section~\ref{Sec:strebel} we review the
notion of Strebel differentials
for smooth curves and give, without proof, two
examples of its degeneration as a curve tends to
a singular stable curve in $\MM_{g,n}$.

In Section~\ref{Sec:stable} we review the notion
of dualizing sheaf for a stable curve and define Strebel
differentials for nonsmooth stable curves. We
show that Strebel differentials with given
perimeters $p_1, \dots, p_n$ form a continuous
section of the vector bundle of quadratic
differentials over $\MM_{g,n}$. This theorem
does not, as far as we know, appear in the
literature. 

In Section~\ref{Sec:compactifications} we
describe the compactification $K\MM_{g,n}$
and the cell complex $A$ homeomorphic to
$K\MM_{g,n} \times \R_+^n$. We prove
that they are indeed homeomorphic, and
also briefly sketch Looijenga's proof~\cite{Looijenga95}.

In Section~\ref{Sec:Chern} we describe a
framework that allows one to work with
differential forms on cell complexes. The
forms representing the connections on the bundles $\B_i$ 
and their curvatures fit into this
framework and thus their use is justified.

\subsection{Acknowledgements}

I am greatful to E.~Looijenga for pointing out the 
references~\cite{Looijenga95},
\cite{Keel}, and~\cite{Keel2}, to M.~Kontsevich,
S.~Lando, C.~Maclean, G.~Shabat, and K.~Strebel for
useful discussions and for their interest in the subject, 
and to A.~Zvonkin for his remarks on the text.

\section{Strebel's theorem} \label{Sec:strebel}

\subsection{The case of smooth curves: a review}
\label{Ssec:smooth}

Let $C$ be a Riemann surface. A {\em simple
differential} on $C$ is a meromorphic
section of its cotangent bundle. In a local coordinate $z$, it
can be written as $f(z) dz$, where $f$ is a meromorphic
function. A {\em quadratic differential}
is a meromorphic section of the tensor square of 
the cotangent bundle. In a local coordinate $z$,
it can be written as $f(z) dz^2$. 

Let $\varphi$ be
a quadratic differential and $z_0 \in C$ a
point that is neither a pole nor a zero of~$\varphi$.
Then $\varphi$ has a square root in the neighborhood
of $z_0$: it is a simple differential $\gamma$,
unique up to a sign, such that $\gamma^2 = \varphi$.
The integral
$$
\int_{z_0}^z \gamma
$$
is a biholomorphic mapping from a neighborhood
of $z_0$ in $C$ to a neighborhood of $0$ in $\mathbb{C}$.
The preimages of the horizontal (vertical) lines in
$\mathbb{C}$ under this mapping are called {\em horizontal
{\rm(}vertical{\rm)} trajectories} of the quadratic differential
$\varphi$. Because $\gamma$ is defined up to a sign
these trajectories do not have a natural orientation.

If $z_0$ is a double pole, then in the neighborhood
of $z_0$ the quadratic differential $\varphi$ has an expansion
$$
\varphi = a \frac{dz^2}{(z-z_0)^2} + \dots \, .
$$
The complex number $a$ is called the {\em residue}
of the double pole and does not depend on the
choice of the local coordinate.

By a local analysis, it is easy to see that if $z_0$ is
a $d$-tuple zero of $\varphi$, then there are $d+2$
horizontal trajectories issuing from $z_0$. Further,
if $z_0$ is a simple pole, then there is a unique
horizontal trajectory issuing from $z_0$. Finally,
if $z_0$ is a double pole whose residue is a negative
real number $a = -(p/2\pi)^2$, then $z_0$ is surrounded by closed
horizontal trajectories (see Figure~\ref{1}~(b)).
These trajectories, together with $z_0$, 
form a topological open disc in $C$.
In the metric $|\varphi|$ all these trajectories
have the same length $p$. The other
possible cases are the case of a double pole
with a positive real or a nonreal residue and 
the case of poles of order greater than~$2$. We 
will not need them, so we leave them to the reader.

Now we are ready to formulate Strebel's theorem.

Let $C$ be a connected (not necessarily compact)
Riemann surface without boundary with $n \geq 1$ distinct
marked and numbered points $z_1, \dots, z_n$. 
We will be interested only in the case where $C$ is
a surface of genus $g$ with a finite number of
punctures. (The punctures and the marked points
are two different things and are pairwise distinct.)
However Theorem~\ref{Thm:Strebel} below is
applicable to any hyperbolic surface $C$, i.e.,
whenever the universal covering of $C \setminus \{ z_1, \dots, z_n \}$
is the Poincar\'e disc.
For example, $C$ can be a ring, a torus with a puncture,
a surface of genus $2$ with a hole, or a surface of
infinite genus. 

In his book on quadratic differentials 
Strebel proves the following theorem
(\cite{Strebel},~Theorem~23.5 and Theorem~23.2 for $n=1$).

\begin{theorem} \label{Thm:Strebel}
For any positive real numbers $p_1, \dots, p_n$ there exists
a unique quadratic differential $\varphi$ on $C$ satisfying the
following conditions. (i)~It has double poles at the marked
points and no other poles. (ii)~The residue at $z_i$
equals $-(p_i/2\pi)^2$. (iii)~If we denote by $D_i$ the disc domain
formed by the closed horizontal trajectories of $\varphi$ surrounding
$z_i$, then 
$$
\bigcup_i \overline{D_i} = C.
$$
\end{theorem}

\begin{definition}
A quadratic differential satisfying the above conditions
is called a {\em Strebel differential.}
\end{definition}

If $C$ is a compact surface of genus $g$, then the nonclosed
horizontal trajectories of a Strebel differential
$\varphi$ form a connected graph embedded into $C$.
All its vertices (situated at the zeroes of $\varphi$)
have degrees $\geq 3$.
Its edges have natural lengths (measured with the length
measure $\sqrt{|\varphi|}$). Its faces are the disc domains
$D_i$ and are in a one-to-one
correspondence with the marked points. The perimeter
of the $i$th face equals $p_i$. Each face $D_i$,
punctured at its marked point, has
a natural flat Riemannian metric $|\varphi|$. In this
metric it is isometric to a semi-infinite cylinder
whose base is a circle of length $p_i$.

\begin{definition}
An {\em embedded}\, or {\em ribbon} graph is a connected
graph endowed with a cyclic order of the half-edges issuing
from each vertex.
\end{definition}

It is a standard fact that any given cyclic order 
allows one to construct a unique embedding of the
graph into a surface. 
If an abstract graph is embedded into a 
surface, the cyclic order of the half-edges adjacent
to a vertex is just the counterclockwise order.

On the set $H$ of the half-edges of the graph we 
introduce three permutations: $\sigma_0$ is the product
of the cyclic permutations assigned to the vertices,
$\sigma_1$ is the involution without fixed points
that exchanges the half-edges of each edge, 
$\sigma_2 = \sigma_0^{-1} \sigma_1$ is the permutation
whose cycles correspond to the faces of the ribbon graph.
These permutations sum up all the information about the ribbon
graph.

\begin{proposition} \label{Prop:strips}
If we are given a ribbon graph with $n$
numbered faces, endowed with edge lengths, and
such that each vertex has a degree at least $3$, then
there is a unique way to recover
a Riemann surface $C$ with $n$ marked points and to determine
the perimeters $p_i$ in such a way that the ribbon graph is
the graph of nonclosed horizontal trajectories 
of the corresponding Strebel differential.
\end{proposition}

The construction given in the proof is described, for
example, in~\cite{Kontsevich}, Section~2.2.

\paragraph{Proof of Proposition~\ref{Prop:strips}.}
To find the perimeters $p_i$ we just add up the lengths of the
edges surrounding each face. 

\begin{figure}[h]
\begin{center}
\
\epsfbox{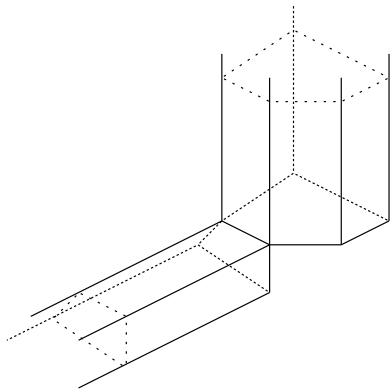}

\caption{Gluing a Riemann surface from strips. 
\label{Fig:strips}}

\end{center}
\end{figure}

The Riemann surface is obtained in the following way.
To every oriented edge of the ribbon graph we assign a strip
$[0,l] \times [0, + i \infty[$ in the complex plane,
where $l$ is the length of the edge. This strip 
inherits the standard complex structure from the complex
plane. Now we construct our surface by gluing together
the strips corresponding to all the oriented edges
(Figure~\ref{Fig:strips}). 

First, for every edge, we glue together the two strips 
that correspond to the two ways of orienting
this edge. The segment $[0,l]$ is identified
with $[l,0]$ and the complex structure is
extended in the natural way. Now we glue together,
along the sides $[0 , + i \infty[$,
the strips that correspond to neighboring edges
in the same face. The complex structure is extended
naturally to $]0, + i \infty[$. It remains to extend
the complex structure to the vertices of the
ribbon graph and to the $n$ punctured points. 

At a vertex of degree $k$ there are $2k$ right angles
of strips that meet together, that is, in whole, an angle
of $k \pi$. Let us place the 
vertex at the origin of the complex plane and put the strips
on the plane one after another around the vertex
(so that the 5th strip will overlap with the
1st one, the 6th one with the 2nd one, and so on). 
If $z$ is the coordinate on the complex plane,
we introduce a local coordinate at the neighborhood 
of the vertex using the function $z^{2/k}$.

Finally, consider a marked point and the semi-infinite
cylinder formed by the strips that surround it. Let
$h$ be the height of a point in this cylinder and
$\theta \in [0, 2\pi[$ its argument (the origin
of the angles can be chosen arbitrarily). Then
$e^{i\theta - h}$ is a local coordinate in the
neighborhood of the marked point.

The uniqueness of the Riemann surface is proved in the
following way. Consider a Strebel differential on
a Riemann surface. Let us
cut the Riemann surface along the nonclosed horizontal
trajectories of the Strebel differential and along its
vertical trajectories joining
the marked points to the vertices of the ribbon graph
(the zeroes of the differential). We obtain the set
of strips described above. Therefore our Riemann surface
is necessarily glued of strips as in the above
construction. \qed

\bigskip

Denote by $B_i$ the polygon (the face of the graph) that forms 
the boundary of the $i$th disc domain $D_i$. (If $D_i$ is
adjacent to both sides of an edge $e$, this edge appears
twice in the polygon $B_i$.)
Further, denote by $T_i$ the complex line tangent to $C$
at $z_i$ and by
$ST_i = (T_i \setminus \{0 \})/\R_+$ its real spherization.
(Here and below $\R_+$ is the set of positive real numbers.)
Then there exists a canonical identification
$$
B_i = ST_i.
$$
Indeed, given a direction $u \in ST_i$,
there is a unique vertical trajectory of $\varphi$ issued
from $z_i$ in the direction $u$. This trajectory
meets the polygon $B_i$ at a unique point, and this point
will be identified with $u$.

Thus Strebel's theorem allows us to define $n$ polygonal
bundles $\B_i$ over $\M_{g,n} \times \R_+^n$ and these bundles can be
identified with the circle bundles obtained by a real
spherization of the complex
line bundles $\L_i^*$. (Recall that the fiber of the bundle
$\L_i$ is the cotangent line to $C$ at $z_i$.)

We are going to show that these polygonal bundles can be
extended to $\MM_{g,n} \times \R_+^n$,
where $\MM_{g,n}$ is the Deligne-Mumford compactification,
so that the identification above is preserved.

\subsection{Two examples}
\label{Ssec:examples}

To extend the bundles $\B_i$ over $\MM_{g,n}$, 
we need to extend to stable curves the notion
of Strebel differentials. 
The construction is carried out in the next section.
Here we just give two examples, without any proofs.

\begin{example}
Consider the case of a torus with one marked point
that degenerates into a sphere with one marked point
and two identified points (see Figure~\ref{FigTorus}; 
the marked point is represented as a black dot). 
Fix a positive real number $p$.
On every torus there exists a unique
Strebel differential with residue $-(p/2\pi)^2$ at the marked point.
It determines a $1$-faced embedded (ribbon) graph composed of the
nonclosed horizontal trajectories. This graph is either
a hexagon whose opposite edges are glued together in pairs, or a
quadrilateral whose opposite edges are glued together in pairs. In the
figure we represented a hexagon.

Now, when the torus degenerates into a sphere, the lengths
of the edges $l_1$ and $l_2$ tend to $0$. Thus, on the
sphere we obtain a graph with only one edge of length
$l_3 = p/2$. This edge joins the two identified points.
If we put the identified points at $0$ and $\infty$,
and the marked point at $1$, then the limit Strebel differential 
on the sphere equals
$$
\varphi = -(p/2\pi)^2 \, \frac{dz^2}{z(z-1)^2}.
$$
It has simple poles at the identified points $0$ and $\infty$.

\begin{figure}
\begin{center}
\
\epsfbox{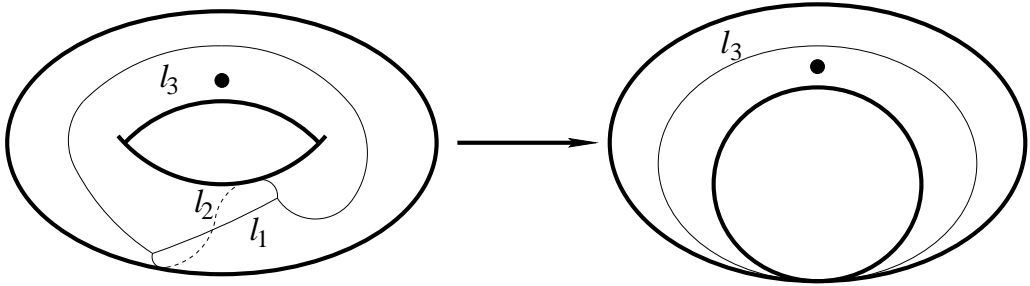}
\end{center}
\caption{A torus degenerating into a sphere with two identified
points. \label{FigTorus}}
\end{figure}

\end{example}

\begin{example}
Now consider a sphere with four marked points that degenerates
into a reducible curve consisting of two spherical components
intersecting at one point (Figure~\ref{FigSphere}). 
Assume that the first component contains
the marked points $z_1$ and $z_2$, while the second component
contains the marked points $z_3$ and $z_4$. Fix 4 positive
real numbers $p_1, p_2, p_3, p_4$. We will assume that 
$p_1 > p_2$, but $p_3 = p_4$ (in order to obtain two different
pictures on the two components).

For any positions of the fixed points on the sphere $\CP^1$,
there is a unique Strebel differential with residues
$-(p_i/2\pi)^2$ at the marked points $z_i$. This differential determines
a $4$-faced graph on $\CP^1$. As the curve approaches the
degeneration described above, this graph necessarily becomes of
a particular form. Namely, it will contain a simple cycle, formed
by several edges, separating the marked points $z_1$ and
$z_2$ from the marked points $z_3$ and $z_4$. In other words,
the faces number $1$ and $2$ become adjacent as well as the
faces $3$ and $4$. When the curve degenerates, the lengths
of all the edges in the above cycle tend to $0$.

On the first component we obtain a graph with 2 vertices.
The vertex at the nodal point has degree $1$,
the second vertex has degree $3$. The corresponding
quadratic differential has a simple pole at the nodal point
and a simple zero at the other vertex (and, of course, double
poles at the marked points).

On the second component we obtain a graph with a unique vertex
of degree $2$ at the nodal point. The corresponding quadratic
differential does not have zeros or poles (except the double
poles at the marked points).

\begin{figure}
\begin{center}
\
\epsfbox{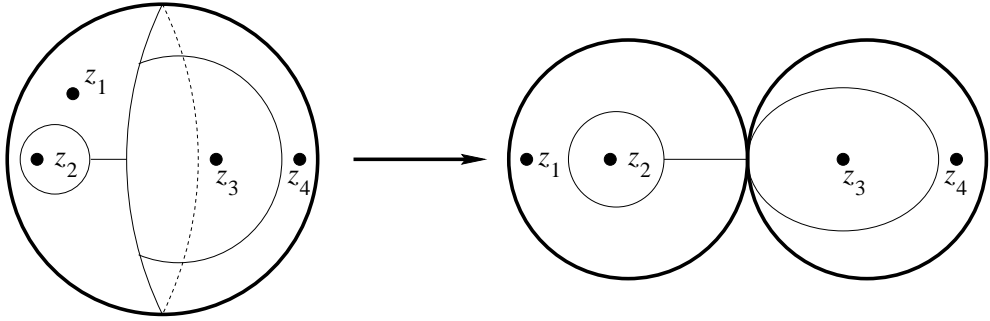}
\end{center}
\caption{A sphere degenerating into a curve with two spherical
components. \label{FigSphere}}
\end{figure}

\end{example}

\section[Differentials on stable curves]
{Simple, quadratic, and Strebel differentials on non\-smooth 
stable curves}
\label{Sec:stable}

\subsection{Simple differentials}
\label{Ssec:simple}

\begin{definition}
Let $C$ be a stable curve with $n$ marked points
$z_i$. A {\em simple differential} $\gamma$ on $C$ is
a meromorphic differential defined on each component of 
$C$ and satisfying the following properties. (i)~It
has at most simple poles at the marked points and at the
nodal points, but no other poles. (ii)~For each nodal point,
the sum of the residues of the poles of $\gamma$ on the
two components meeting at this point vanishes.
\end{definition}

One can readily check that simple differentials form
a vector space $V$ of dimension $g+n-1$ (if $n \geq 1$)
for any stable curve of arithmetic genus $g$, 
whether it is smooth or not. 

Indeed, consider a stable curve with several irreducible
components $C_i$. Suppose $C_i$ is of genus
$g_i$, has $n_i$ marked points and $m_i$ nodal points.
Then we have
$$
n = \sum n_i, \qquad 2-2g =  \sum (2 - 2g_i - m_i)
$$

Now, according to the
Riemann-Roch theorem, the dimension of the space of
sections of a line bundle with first Chern class $c$, such
that at most simple poles are allowed at $k$ fixed points,
is equal to $c+k+1-g$, whenever $c+k \geq 2g-1$.
In our case, on $C_i$, $c = 2g_i-2$ (the
first Chern class of the cotangent line bundle)
and $k = n_i + m_i$. Thus the dimension of the 
space of sections equals $g_i+n_i+m_i-1$.
Adding these numbers for all the irreducible components
$C_i$ and subtracting the total number 
$\frac12 \sum m_i$ of nodal points
(because each nodal point gives a linear relation on the
residues) we obtain $g+n-1$.

Because the dimensions of the spaces $V$ are the same,
these spaces form a holomorphic
vector bundle $\ocV$ over the space $\MM_{g,n}$.
This follows immediately from algebro-geometric arguments.
Indeed, denote by ${\overline {\cal C}}_{g,n}$ the universal curve
over $\MM_{g,n}$. Then simple differentials form
a sheaf on ${\overline {\cal C}}_{g,n}$. It is
called the {\em relative dualizing sheaf} and is the
sheaf of sections of a line bundle (the relative
cotangent line bundle) over ${\overline {\cal C}}_{g,n}$.
(See, for example,~\cite{Hartshorne}, Chapter~III, Theorem~7.11,
where it is proved that the dualizing sheaf is the sheaf
of sections of a line bundle for any algebraic variety
that is locally a complete intersection, and an explicit
construction of the sheaf is given. See also~\cite{HarMor},
Chapter~3, Section~A.)
The direct image of the relative dualizing sheaf
on $\MM_{g,n}$ has the property that the dimensions
of its fibers are the same. Therefore the direct image
is itself the sheaf of sections of a vector bundle
(see~\cite{Hartshorne}, Exercise 5.8).
 
The fact that the spaces $V$ form a holomorphic vector
bundle can also be understood more intuitively.
The space of pairs $(C, \gamma)$, where $C$
is a stable curve and $\gamma$ a simple differential on it,
can be given a complex structure
in the following way. For any pair $(C,\gamma)$,
one can calculate the integral of $\gamma$ over any closed
loop that does not contain marked or nodal points.
The complex structure is introduced
by requiring that these integrals be meromorphic functions on the
space of pairs $(C,\gamma)$. (That these integrals
can have poles is shown by the following example.
Consider a torus whose meridian is contracted, so that it
degenerates into a sphere with two identified points. Then
the integral of a simple differential over the parallel
will tend to $\infty$.)

\subsection{Quadratic differentials}
\label{Ssec:quadratic}

Now we repeat the above construction for quadratic differentials.

\begin{definition}
Let $C$ be a stable curve with $n$ marked points
$z_i$. A {\em quadratic differential} $\varphi$ on $C$ is
a meromorphic quadratic differential defined on each component of 
$C$ and satisfying the following properties. (i)~It
has at most double poles at the marked points and at the
nodal points, but no other poles. (ii)~For each nodal point,
the residues of the poles of $\varphi$ on the
two components meeting at this point are equal.
\end{definition}

\begin{remark}
The residue of a quadratic differential at a
pole of order at most $2$ is equal to the coefficient
of $dz^2/z^2$ (for any local coordinate $z$). If the
order of the pole is actually less than $2$, we let the
residue be equal to $0$.
\end{remark}

As above, the dimension of the space $W$ of quadratic
differentials is the same for any stable curve $C$
and equals $3g-3+2n$. 

Indeed, if $C_i$ is an irreducible component
of $C$, that has genus $g_i$ and contains
$n_i$ marked points and $m_i$ nodal points,
then the dimension of the space of quadratic differentials
on it is $3g_i-3 + 2n_i + 2 m_i$. Adding these numbers
for all the components and subtracting the total number
$\frac12 \sum m_i$ of the nodal points (because each nodal
point gives a linear relation on the residues), we obtain
$3g-3+2n$. 

Since the dimensions of the spaces $W$ are the same, they
form a holomorphic vector bundle $\ocW$ over $\MM_{g,n}$. 
This follows from the same arguments as for simple
differentials. The quadratic differentials form a sheaf
on the universal curve ${\overline {\cal C}}_{g,n}$:
the sheaf of sections of the tensor square of the
relative cotangent line bundle. The direct image of
this sheaf on $\MM_{g,n}$ has the property that all
its fibers are of the same dimension. Therefore it
is a sheaf of sections of a holomorphic vector bundle.

\subsection{Strebel differentials}
\label{Ssec:stabStreb}

Here we define Strebel differentials on stable
curves.

Let $C$ be a stable curve with $n$ marked points.
Suppose we are given $n$ positive real numbers $p_1, \dots, p_n$.

\begin{definition} \label{Def:unmarked}
We say that an irreducible component of a stable curve $C$
is {\em marked} if it contains at least one
marked point and {\em unmarked} if it contains no
marked points.
\end{definition}

\begin{definition} \label{Def:StabStreb}
A {\em Strebel differential $\varphi$ on a stable curve $C$}
is a quadratic differential on $C$ satisfying the following properties.
(i)~It has double poles at the marked points, at most
simple poles at the nodal points, and no other poles.
(ii)~The residue of the pole at the $i$th marked point
$z_i$ equals $-(p_i/2\pi)^2$. (iii)~The differential $\varphi$
vanishes identically on the unmarked components. 
(iv)~Let $C'$ be a marked component of
$C$. Let us puncture $C'$ at the nodal points. 
For $z_i \in C'$, denote by $D_i$ the
disc domain formed by the closed horizontal trajectories 
of $\varphi$ surrounding $z_i$. Then we have
$$
C' = \bigcup_{i|z_i \in S} \overline{D_i}.
$$
\end{definition}

\begin{remark}
Strebel differentials have at most simple poles at the
nodes of $C$ (unlike generic quadratic differentials,
that have double poles). Therefore the condition
that the residues of the poles on the
two components meeting at a nodal point must be
equal is automatically satisfied, since both
residues vanish.
\end{remark}

\begin{remark} \label{Rem:stabStreb}
It follows from Strebel's theorem, that (once the positive
numbers $p_1, \dots, p_n$ are given) there exists a unique
Strebel differential $\varphi$ on any stable curve $C$. 
Its restriction to unmarked components vanishes. 
Its restriction to each marked component $C'$ is the 
Strebel differential on $C'$ with punctures at 
the nodal points. Indeed, it is easy to see
that when we put the punctures back into the component $C'$,
the corresponding Strebel differential will have at most
simple poles at these points, because there is only a
finite number of nonclosed horizontal trajectories
issuing from them.
\end{remark}

\begin{remark}
Let $C_1$ be a smooth compact Riemann surface with $n$
marked points, and let $C_2$ be obtained by puncturing
$C_1$ at a finite number of points (different from the
marked points). Given a list of positive real parameters
$p_1, \dots, p_n$, there is a unique Strebel differential
$\varphi_1$ on $C_1$ and a unique Strebel
differential $\varphi_2$ on $C_2$. At first sight, one could
think that they are the same; but, in general, this is not true. Indeed,
if we restrict $\varphi_1$ to $C_2$, we will see that the disc
domains of $\varphi_1$ contain punctures at their interior,
which is not allowed for a Strebel differential. Conversely, if
we try to extend $\varphi_2$ to $C_1$ by putting back the
punctured points, we will, in general, obtain
simple poles at these points. Again, a Strebel differential
is not allowed to have poles outside the marked points.

Thus in the condition (iv) in
Definition~\ref{Def:StabStreb} above,
it is important to puncture each component $C'$ at the
nodal points.
\end{remark}

As in the case of smooth compact Riemann surfaces, the nonclosed
horizontal trajectories of a Strebel differential
on a stable curve form a graph, embedded into the
stable curve. More precisely, it is embedded
into the union of the marked components of
the stable curve. The vertices of the graph are of degrees $\geq 3$,
except for the vertices that lie at the nodal points
and can have any degree $\geq 1$. Its edges have natural
lengths (measured, as before, with $\sqrt{|\varphi|}$). Its faces
are in a one-to-one correspondence with the marked
points, and the perimeter of the $i$th face
equals $p_i$. As before, if we denote by $B_i$ the polygon
surrounding the $i$th face
of the graph, by $T_i$ the complex line tangent to the
marked point $z_i$, and
by $ST_i = (T_i \setminus \{ 0\})/\R_+$ its real spherization,
we have a canonical identification
$$
B_i = ST_i.
$$

\subsection{Stable ribbon graphs}
\label{Ssec:stabribgraphs}

In Section~\ref{Sec:compactifications} 
we will need a formal definition
of a graph formed by the nonclosed horizontal trajectories
of a Strebel differential on a stable curve. We give
the definition here. 

To understand the definition below
one must imagine that we have contracted to a point
each unmarked component of the stable curve. Thus we have 
obtained a graph embedded into a new (usually singular) curve. 

The unmarked components of the initial
stable curve form a not necessarily connected
subcurve in it. Each connected component of this subcurve is 
contracted to a vertex of our graph.
On each such vertex we will mark the
arithmetic genus of the corresponding contracted 
component. 

If a vertex $v$ of the graph lies at a node of the curve,
there is no longer any natural cyclic order on the
half-edges issuing from it. Instead, we will have a permutation
(with several cycles) acting on these half-edges. Each
cycle of this permutation corresponds to a component
of the curve at the neighborhood of $v$. The cycle
determines the counterclockwise cyclic order of the half-edges on the
corresponding component.

\begin{definition} \label{Def:stabgraph}
A {\em stable ribbon graph} is a connected
graph endowed with the following structure. 
(i)~A non-negative integer (a {\em genus defect})
is assigned to each vertex.
(ii)~A permutation acting on the set of
half-vertices issuing from each vertex is given.

There are two types of vertices whose genus defect cannot
be equal to~$0$: first, the vertices of degree~$1$, second,
the vertices of degree two such that the corresponding
permutation of the two half-edges is a transposition.
\end{definition}

\begin{figure}[h]
\begin{center}
\
\epsfbox{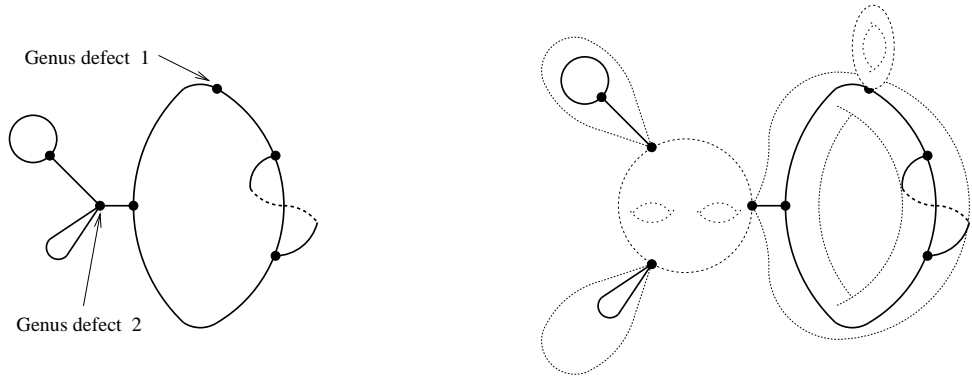}

\caption{A stable ribbon graph and its
surface of embedding. \label{Fig:stabgraph}}
\end{center}
\end{figure}

Kontsevich (\cite{Kontsevich}, Appendix B) 
gives an equivalent definition of stable ribbon graphs.
A stable ribbon graph is represented in 
Figure~\ref{Fig:stabgraph}. Its surface of embedding
(that can be uniquely reconstructed from the 
stable ribbon graph structure) is shown in dotted lines.

It is easy to define faces of a
stable graph. Let $H$ be the set of all the half-edges
of the graph, $\sigma_0$ the product of all the permutations
(with disjoint supports) assigned to the vertices,
and $\sigma_1$ the involution without fixed
points that switches two half-edges of each edge.
Then a face is a cycle of the permutation 
$\sigma_2 = \sigma_0^{-1} \sigma_1$.
The permutations $\sigma_0$, $\sigma_1$, and $\sigma_2$ sum up
all the structure of a stable ribbon graph, except the
genus defect function.
We will usually consider stable graphs with $n$
numbered faces.

The {\em genus} of a stable ribbon graph is the arithmetic
genus of its surface of embedding, plus the sum of genus
defects of all its vertices. 

In Section~\ref{Sec:compactifications} we will see
that stable graphs are obtained from 
ordinary ribbon graphs by edge contractions.

\begin{proposition} \label{Prop:strips2}
Given a stable ribbon graph with $n$ numbered faces and
endowed with edge lengths, we can find a set
of perimeters $p_1, \dots, p_n$ and a stable curve
such that the stable ribbon graph is the graph of
nonclosed horizontal trajectories of the corresponding
Strebel differential. The marked components of the
curve are uniquely determined.
\end{proposition}

\paragraph{Proof.} To find the perimeters we simply
add up the lengths of the edges surrounding each face.
To construct a stable curve we carry out the same operations
of strip gluings as in the proof of Proposition~\ref{Prop:strips}.
This gives us the complex structure on the marked components.
That it is unique is shown as in Proposition~\ref{Prop:strips}.
As for the unmarked components, their arithmetic genera are
given by the genus defect function of the stable ribbon
graph, but their complex structures and even
their topologies can be chosen arbitrarily. \qed

\subsection{Strebel differentials form a continuous family}
\label{Ssec:smoothfamily}

Here we prove the continuity of the map that assigns
to a stable curve and a list of perimeters in
$\MM_{g,n} \times \R_+^n$ the corresponding
Strebel differential in the total space of the vector 
bundle $\ocW$ of quadratic differentials over $\MM_{g,n}$. 
In particular, Strebel differentials with fixed perimeters form
a continuous section of the vector bundle $\ocW$ over $\MM_{g,n}$. 
This fact does not seem to be
stated explicitly in the literature, although it would not
be a surprise for a specialist in Teichm\"uller spaces.

We will need a rather well-known characterization
of convergence in $\MM_{g,n}$ and in the total
space of the vector bundle $\ocW$ of quadratic differentials.

We view a complex structure on a surface as an operator
$J$ that multiplies each tangent vector by $i$.
We also remind the reader that quadratic differentials have at most
double poles at the marked points, at most double poles with equal
residues at the nodes, and no other poles.

\begin{definition} \label{Def:deformation}
A continuous map $f:C_1 \rightarrow C_2$ from
a stable curve to another one is called a {\em deformation}
if (i)~the preimage of any node of $C_2$ is either a node
of $C_1$ or a simple loop in the smooth
part of $C_1$, (ii)~$f$ is an orientation-preserving
diffeomorphism outside the nodes and loops,
and (iii)~$f$ sends the marked points to the
marked points preserving their numbers.
\end{definition}

For a summary of various properties of deformations
and their relations with augmented Teichm\"uller spaces
see~\cite{Bers}. See also~\cite{Abikoff} for related
questions on the topology of the Teichm\"uller spaces.

\begin{proposition} \label{Prop:teichmetric}
Let $(C_m, \varphi_m) \rightarrow (C, \varphi)$ 
be a converging sequence in $\ocW$ and denote by $J_m$ and $J$
the complex structures on $C_m$ and $C$. 
Starting from some $m$ there exists 
a sequence of deformations $f_m:C_m \rightarrow C$
such that:

(i)~on any compact set $K \subset (C \setminus \mbox{\rm nodes})$ 
the sequence of complex structures $(f_m)_* J_m$ 
converges uniformly to $J$;

(ii)~on any compact set 
$K \subset (C \setminus \mbox{\rm nodes and marked points})$ 
the sequence of complex-valued symmetric $2$-forms $(f_m)_* \varphi_m$
converges uniformly to $\varphi$.
\end{proposition}

\paragraph{Sketch of a proof.}
Consider a point $x \in \MM_{g,n}$ and the corresponding
stable curve with marked points $C_x$. Let $(U, \Stab \, x)$
be a sufficiently small chart containing $x$, where $U$ is an open ball
in $\C^{3g-3+n}$ and $\Stab \, x$ a finite group acting on $U$
and stabilizing $x$. 
The neighborhood of $x$ in $\MM_{g,n}$ is identified
with $U/\Stab \, x$. Consider the part of the
universal curve ${\ocC}_U$ that lies
over $U$. It is a fiber bundle over $U$ whose fibers
are stable curves parameterized by the points of $U$. 
Even as a smooth manifold $\ocC_U$ is, of course,
not a direct product with $U$,
because its fibers have different topologies. However, there exists a
continuous function $f: \ocC_U \rightarrow C_x$ with the
following properties. 
(i)~The restriction of $f$ to any fiber $C_y$ is a deformation
$C_y \rightarrow C_x$ in the sense of 
Definition~\ref{Def:deformation}.
(ii)~The function $f$ is smooth on $\ocC_U$ outside the preimages
of the nodes of $C_x$.
(iii)~The restriction of $f$ to $C_x$ is the identity map.

The function $f$ is a kind of universal family of deformations
over the open set $U$. Such a function can be constructed, 
for example, in the following way.
First of all, let us choose the loops to be contracted by $f$
in each fiber $C_y$. Their free homotopy types 
are uniquely determined by the property that
we must obtain the curve $C_x$ by pinching all these
loops. We choose the loops themselves to be the shortest
geodesics inside the corresponding homotopy classes, with respect
to the unique complete metric of curvature $-1$ on $C_y$, 
compatible with the conformal structure. 
Erasing all the loops and the nodes in each fiber 
we obtain a locally trivial
fiber bundle over a contractible base $U$. Therefore we can
trivialize it by a diffeomorphism
$$
(\ocC_U \setminus \mbox{nodes and loops})
\; \rightarrow \; U \times (C_x \setminus \mbox{nodes})
$$
commuting with the projections to $U$.
We take $f$ to be the second component of this 
diffeomorphism and we extend it to the loops in every
fiber $C_y$ by sending them to the corresponding nodes 
of $C_x$.

For $m$ big enough, $C_m$ lies in $U$ (or,
more precisely, in $U/\Stab \, x$, but we can choose any lifting
of $C_m$ to $U$). For a compact set
$K \subset (C \setminus \mbox{nodes})$
there is, on the whole set $f^{-1}(K)$, a smooth
linear operator $J_U$ acting on tangent planes to the fibers
of $\ocC_U$. Therefore the sequence $f_* J_m$ converges
uniformly on $K$ to the complex structure $J$ of
the stable curve $C_x$. This proves Assertion~(i) of the
proposition.

Now consider a holomorphic section of the vector bundle
$\ocW_U$ over $U$. It is represented by a holomorphic section
of a line bundle over the universal curve $\ocC_U$, namely, 
of the tensor square of the relative
dualizing bundle. Almost each fiber
of this line bundle is naturally identified with the
tensor square of the cotangent line to the corresponding
stable curve at the corresponding point. The only exceptions
are the fibers over the marked and the nodal points. 
First assume for simplicity that our sequence $(C_m, \varphi_m)$ 
belongs to (or, more precisely, is a restriction of)
some holomorphic section of $\ocW_U$.
Then, exactly as before, we conclude that if 
$K \subset (C_x \setminus \mbox{nodes and marked points})$ is a compact set,
the sequence of quadratic differentials $f_* \varphi_m$ converges
uniformly on $K$ to the quadratic differential $\varphi$.
In general, the sequence $(C_m, \varphi_m)$ does not belong to
a holomorphic section of $\ocW_U$. Then we have to consider a family
of ${\rm rk} \, \ocW$ holomorphic sections of $\ocW_U$ over $U$,
forming a basis of each of its fibers. The coordinates of the
elements of the sequence $(C_m, \varphi_m)$ in the
basis formed by the sections converge to the
coordinates of $(C, \varphi)$. Applying to each
section of the family the above argument, we conclude that sequence
$f_* \varphi_m$ converges to $\varphi$ uniformly on $K$. 
This proves Assertion~(ii) of the proposition. \qed

\begin{theorem} \label{Thm:contfam}
The Strebel differentials with fixed parameters 
$p_1, \dots, p_n$ form a continuous nonvanishing section
of the vector bundle $\ocW$ of quadratic differentials
over the Deligne-Mumford compactification $\MM_{g,n}$.
\end{theorem}

\paragraph{Proof.} Let us fix a stable curve $C$
and consider a sequence of smooth curves $C_m$ that
tends to $C$ in $\MM_{g,n}$ as $m$ tends to $\infty$.
By Strebel's theorem (Theorem~\ref{Thm:Strebel}) and
Remark~\ref{Rem:stabStreb} there is a unique Strebel
differential $\varphi_m$ on each curve
$C_m$ and a unique Strebel differential $\varphi$ on
$C$. We will prove that
$\varphi_m$ tends to $\varphi$ in the vector bundle
$\ocW$, as $m$ tends to $\infty$.
This is enough to prove the theorem, since smooth
curves form an open dense subset of $\MM_{g,n}$.

First of all, the sequence of quadratic differentials
$\varphi_m$ is bounded, therefore it has
at least one limit point $\tphi$. We will prove that
$\tphi = \varphi$, which implies that the limit point
is unique and is therefore a true limit. To do that,
we study the limit quadratic differential $\tphi$
and prove that it has all the properties of a Strebel
differential.

For shortness we will call just ``trajectories'' the 
horizontal trajectories of the differentials.

{\bf 1. The nonclosed trajectories of $\tphi$ have a
finite total length.} 

Let $x \in C$ be
a regular point that is neither a zero nor a pole of
$\tphi$. Let $x_m \in C_m$ be a sequence of
points of the curves $C_m$ that tends to $x$
in the universal curve $\ocC_{g,n}$ over $\MM_{g,n}$. By moving
each $x_m$ slightly inside $C_m$ we can assume
that each $x_m$ belongs to a closed horizontal trajectory
of $\varphi_m$, because the union of closed trajectories
is dense in each curve $C_m$. Moreover, by extracting
a subsequence we can assume that all these closed
trajectories belong to the disc domains $D_i$ for
the same $i$. Therefore each closed trajectory has
the same length $p_i$. We will prove that $x$ is
contained either in a closed trajectory
or in a nonclosed trajectory shorter than $p_i$.

Suppose that moving along the trajectory through $x \in S$ 
we have covered a segment of length $l > p_i$ without
encountering a nonregular point (a pole, a node, or a zero) 
and without passing
twice through the same point. The above segment has a
compact neighborhood $K$ that does not contain marked points,
nodes of $C$, or zeroes of $\tphi$. Using 
Proposition~\ref{Prop:teichmetric} we can construct a
sequence of deformations $f_m:C_m \rightarrow C$
in such a way that the sequences of complex structures 
$(f_m)_* J_m$ and of quadratic differentials 
$(f_m)_* \varphi_m$ converge uniformly on $K$. Therefore, 
for $m$ big enough, the trajectory through $x_m$ of the
quadratic differential $\varphi_m$ will also have a segment
of length greater than $p_i$. This is a contradiction.
Thus, if $x$ is a regular point of $\tphi$, the trajectory through
$x$ is either closed or nonclosed of finite length. 

By compactness of $C$, a  nonclosed
trajectory of finite length necessarily has two endpoints
in $C$. These endpoints can be zeroes of $\tphi$, simple
poles of $\tphi$ (including possible simple poles at the
nodes of $C$), or nodes of $C$ at which $\tphi$ has 
no poles. Since the number of such points
is finite and there is only a finite number of nonclosed
trajectories issuing from each of them, it follows that the
total number of nonclosed trajectories is finite.
 
{\bf 2. Each closed trajectory of $\tphi$ bounds a disc with
a unique marked point inside it.}

If the trajectory $\alpha$
through $x$ is closed, consider a compact tubular
neighborhood $K$ of $\alpha$ that does not contain marked
points, nodes of $C$, and zeroes of $\tphi$. 
We construct a sequence of deformations $f_m:C_m \rightarrow C$ as in
Proposition~\ref{Prop:teichmetric}. 

For $m$ big enough $f_m^{-1}(K)$ 
contains a closed trajectory $\alpha_m$ of the Strebel differential 
$\varphi_m$. Indeed, let $l$ be a real number greater
than any of the perimeters $p_i$. If we choose $m$ big 
enough and a point $x_m \in C_m$ close enough to $f_m^{-1}(x)$,
a segment of length $l$ of the trajectory of $\varphi_m$ through $x_m$
will be entirely contained in $f_m^{-1}(K)$. But for a generic
choice of $x_m$ the trajectory through $x_m$ is closed of
length less than $l$. Therefore we have obtained a closed
trajectory $\alpha_m$ entirely contained in $f_m^{-1}(K)$. Its homotopy
type is uniquely determined by $K$. Indeed, $\alpha_m$
can neither bound a small disc inside $f_m^{-1}(K)$ (because
$\varphi_m$ has no poles inside $f_n^{-1}(K)$), nor have
self-intersections.

We know that this closed trajectory belongs to a disc domain
$D$ of $\varphi_m$ and that the restriction of $f_m$ to $D$ is a
diffeomorphism. Thus $\alpha$, just as $\alpha_m$, 
surrounds a disc that contains a unique marked point.

Note that we have always assumed that $x$ is not a zero
of $\tphi$. Therefore we still know nothing about the
existence of irreducible components of $C$ on which $\tphi$
vanishes identically.

{\bf 3. The poles of $\tphi$ at the nodes of $C$
are at most simple.} 

\begin{figure}
\begin{center}
\
\epsfbox{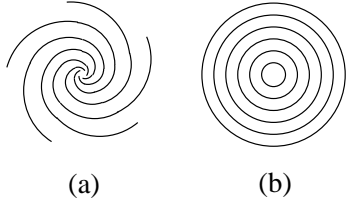}
\end{center}
\caption{Horizontal trajectories near a double pole for: (a) not 
a real negative residue; (b) a real negative residue. \label{1}}
\end{figure}

Consider a nodal point $x$ of $C$, and let us prove that
the pole of $\tphi$ at $x$ is not double, as for a
generic quadratic differential, but at most simple (on both
components meeting at $x$). Suppose that this is not true,
and both poles of $\tphi$ are double (and have the same
residue, as is always the case for quadratic differentials).
Then the common residue is necessarily a negative
real number. Indeed, otherwise a neighborhood of $x$ 
entirely consists of nonclosed trajectories of $\tphi$ of
infinite lengths (see Figure~\ref{1}a), which
contradicts {\bf 1}. If the common residue
is a negative real number, then the point $x$ is surrounded
(on both components) by concentric closed
trajectories (Figure~\ref{1}b). Each of these trajectories
must surround a unique marked point. But this
would mean that $C$ is composed of $2$ 
spherical irreducible
components with one marked point and one nodal point
of each. This is impossible because the curve $C$
is stable.  Thus $\tphi$ cannot have double poles at nodal points.

\begin{figure}
\begin{center}
\
\epsfbox{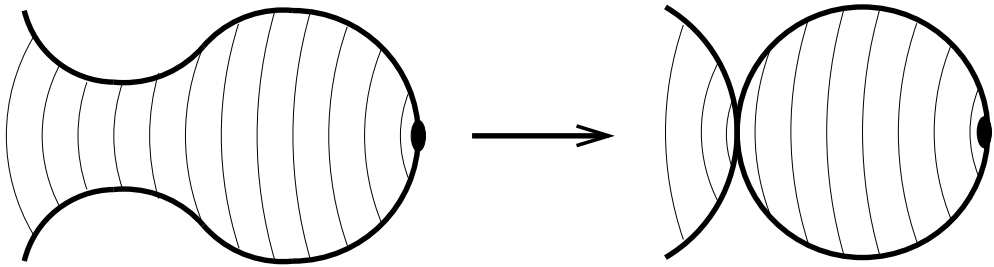}
\end{center}
\caption{Double poles at a nodal point cannot arise
from a disc domain because it would mean that the 
curve $C$ is not stable. \label{2}}
\end{figure}

\begin{figure}
\begin{center}
\
\epsfbox{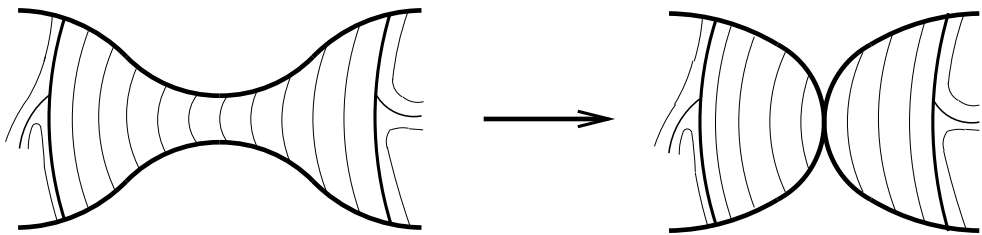}
\end{center}
\caption{Double poles at a nodal point cannot arise
from a cylindric domain because the Strebel differentials
$\varphi_m$ do not have cylindric domains.  \label{3}}
\end{figure}

To get a better insight why double poles are impossible at
the nodes of $C$, we have represented (Figures~\ref{2} and~\ref{3}) 
two families of curves with quadratic differentials, 
degenerating to a nodal curve
$C$ on which the limit quadratic differential has a double
pole at a node. All the quadratic differentials in question have only
a finite number of nonclosed horizontal trajectories. However, in the
first case, the limit curve $C$ is not stable, while in the
second case, the quadratic differentials before the limit
have cylindric domains.

{\bf 4. On components without marked points we have
$\tphi = 0$.} 

Now consider an unmarked component $C'$ of $C$
(see Definition~\ref{Def:unmarked}). Suppose $\tphi$
is not identically equal to zero on this component.
Then $\tphi$ has only a finite number of nonclosed
trajectories on $C'$ and these are of finite
lengths. But, on the other hand, $\tphi$ has no
closed trajectories, because each 
closed trajectory surrounds a marked point and $C'$ 
contains no marked points. This is a contradiction.

{\bf 5. The marked components are covered
by the closures of the disc domains of $\tphi$.}

Finally, consider the restriction of $\tphi$ to a marked
component $C'$ (see Definition~\ref{Def:unmarked}). 
It has double poles with
residues $-(p_i/2\pi)^2$ at the marked points 
(and therefore it does not vanish). 
It has at most simple poles at the nodal points. Each of its
closed trajectories surrounds a unique marked point and
therefore belongs to a disc domain.
The total length of its nonclosed trajectories
is finite and therefore $C'$ is covered
by the closures of the disc domains. 
Thus (by Strebel's theorem)
$\tphi$ is the unique Strebel differential on 
$(C' \setminus \mbox{nodes})$ with parameters $p_i$.

We have proved that $\tphi = \varphi$.

This completes the proof. \qed

\begin{theorem} \label{Thm:cont2}
The map $\MM_{g,n} \times \R_+^n \rightarrow \ocW$ that
assigns to a stable curve and a list
of perimeters the corresponding Strebel differential
is continuous.
\end{theorem}

\paragraph{Proof.} We consider a sequence of
smooth curves $C_m$ tending to a stable
curve $C$, together with a sequence of $n$-tuples
of positive real numbers
$(p_1^{(m)}, \dots, p_n^{(m)})$ tending to
an $n$-tuple $(p_1, \dots, p_n)$ of positive
real numbers. Now we repeat the proof of
Theorem~\ref{Thm:contfam} without modifications.
\qed

\begin{remark} K.~Strebel (\cite{Strebel}, Theorem~23.3)
proves that Strebel differentials on any connected, 
not necessarily compact Riemann surface 
(of finite type and without boundary) with $n$ marked points
depend continuously on the parameters $p_1, \dots, p_n$,
in the topology of uniform convergence on compact sets
outside the marked points.
\end{remark}

\section{The ``minimal reasonable compactification'' of $\M_{g,n}$}
\label{Sec:compactifications}

Here we describe a compactification of $\M_{g,n}$
that is different
from the Deligne-Mumford one. This compactification,
multiplied by $\R_+^n$, is
isomorphic (as a topological orbifold) to a natural
closure of the cell decomposition of $\M_{g,n} \times \R_+^n$
given by Strebel differentials.

\subsection{Cell complexes}

Strebel's theorem allows one to divide the space 
$\M_{g,n} \times \R_+^n$ into cells: two Riemann surfaces
endowed with perimeters $(C; p_1, \dots, p_n)$ and
$(C'; p_1', \dots, p_n')$ belong to the same cell if
the nonclosed horizontal trajectories of
the corresponding Strebel differentials form isomorphic
ribbon graphs (without taking into account the
lengths of the edges). The cell corresponding to a
ribbon graph $G$ whose set of edges is $E$, 
is isomorphic to $\R_+^E/\Aut(G)$,
where $\Aut(G)$ is the (finite) group of automorphisms of the
ribbon graph. A cell $\Delta_1$ is a face of another cell
$\Delta_2$ iff the corresponding graph $G_1$ can be obtained from the graph
$G_2$ by contracting several edges. Gluing such cells
together we obtain an orbifold cell complex
homeomorphic to $\M_{g,n} \times \R_+^n$.

Now we construct a bigger cell complex
whose cells correspond to stable ribbon graphs of genus $g$
with $n$ numbered faces.
Let us first define the operation of edge contracting
in stable ribbon graphs. A cell $\Delta_1$ of our new 
cell complex will be a face of another cell
$\Delta_2$ iff the corresponding stable graph $G_1$ 
can be obtained from the stable ribbon graph
$G_2$ by contracting several edges. 

All the graphs considered below are stable ribbon graphs
of genus $g$ with $n$ numbered faces 
(see Definition~\ref{Def:stabgraph}). Ordinary ribbon graphs with $n$
numbered faces and such that their vertices have degrees at
least~$3$ are particular cases of stable ribbon graphs
(the genus defect at all vertices being equal to~$0$).

Let $G$ be a stable ribbon graph and $e$ its edge. We suppose
that $e$ does not constitute a face on its own. Recall that
$\sigma_0$ is the permutation of the half-edges obtained
by multiplying all the permutations assigned to
vertices; $\sigma_1$ is the involution exchanging
the half-edges of each edge;
$\sigma_2 = \sigma_0^{-1} \sigma_1$ is the permutation
whose cycles correspond to faces.

\begin{figure}[h]
\begin{center}
\
\epsfbox{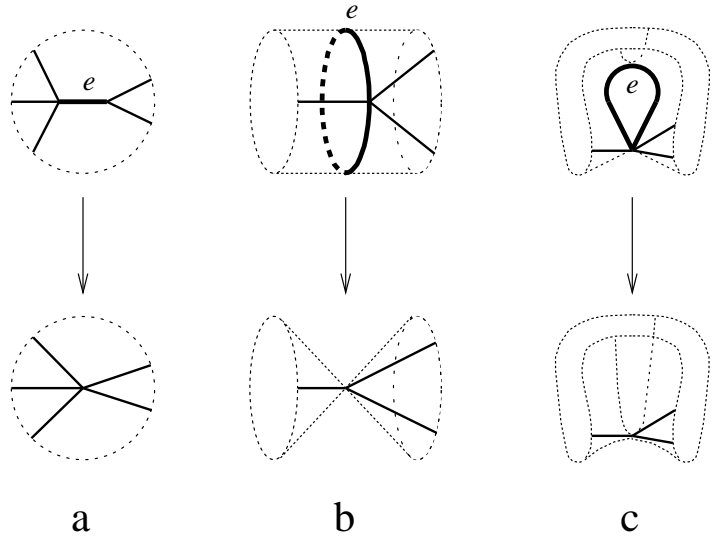}

\caption{Contracting an edge $e$ of $G$. We have represented
the neighborhood of the edge $e$ in the involved
component of the surface of embedding of $G$.
\label{Fig:contraction}}
\end{center}
\end{figure}

\begin{definition}
The {\em contraction of the edge}\, $e$ in the stable
ribbon graph $G$ gives the following stable ribbon graph. 

The underlying combinatorial graph is just the underlying
graph of $G$ with the edge $e$
contracted. 

If the edge $e$ is not a loop (Figure~\ref{Fig:contraction}a), 
then the genus defect assigned to the vertex obtained by its 
contraction is the sum of the genus defects of the
two initial vertices of $e$. If $e$ is a loop based
at a vertex $v$ and the two
half-edges of $e$ belong to the same cycle of the permutation
$\sigma_0$ (Figure~\ref{Fig:contraction}b), 
then the genus defect of $v$ does not 
change after the contraction. If $e$ is a loop
and the half-edges of $e$
belong to two different cycles of $\sigma_0$
(Figure~\ref{Fig:contraction}c), 
then the genus defect increases by $1$.

Finally, the new permutations $\sigma_0'$, $\sigma_1'$,
$\sigma_2'$ are defined as follows. Let $h$ be a half-edge.
Then the half-edge $\sigma_2'(h)$ is the first among
the half-edges $\sigma_2(h)$, $\sigma_2^2(h)$, \dots
that is not a half-edge of $e$.
(In other words, from the point of view of a face whose
boundary included the edge $e$, this edge simply got contracted.)
The permutation $\sigma_1'$ is defined in the obvious way
(by excluding from $\sigma_1$ the cycle corresponding to $e$). The
permutation $\sigma_0'$ equals $\sigma_1' \sigma_2'^{-1}$.
\end{definition}

This operation of edge-contracting might look complicated,
but actually it is quite natural and describes what
happens to the graph of nonclosed horizontal trajectories
of a Strebel differential as the length of one of its
edges $e$ tends to~$0$. From the point of view of each polygon
$B_i$ surrounding a face nothing special happens: if $e$
was part of $B_i$ it simply gets contracted. Simple topological
considerations allow one to find what happens to the
genus defect. The precise statement of this is given
in Theorem~\ref{Thm:homeomorphism} below.

Note that if an edge $e$ is the unique
edge that surrounds a face, then it cannot be
contracted, because its length must remain equal
to the perimeter $p_i$ of the corresponding face.

The two propositions below are simple combinatorial
exercises.

\begin{proposition}
Edge contracting is commutative,
in other words the result of a contraction of two edges
does not depend on the order in which they are performed.
\end{proposition}

Thus it makes sense to talk about contracting a subset of
the set of edges of a stable ribbon graph.

Consider a stable ribbon graph $G$. Let $E$ be the set
of its edges, $H$ the set of its half-edges, and
$\sigma_0, \sigma_1, \sigma_2$ its structural permutations. Decompose
$E$ into a disjoint union $E = E_c \sqcup E_r$. We are
going to describe the result of the contraction of the edges
of $E_c$. 

Decompose $E_c$ into the union of connected
components, $E_c = E_1 \sqcup \dots \sqcup E_k$. First
we introduce a structure of
a stable ribbon graph on each of the $E_i$. The underlying
graph is just the subgraph of $G$ with edges $E_i$. The
permutation $\sigma_1^{(i)}$ is the restriction of
$\sigma_1$ to the half-edges of $E_i$. The image of
a half-edge $h$ under $\sigma_0^{(i)}$ is the first
half-edge among $\sigma_0(h)$, $\sigma_0^2(h)$, \dots
to belong to an edge of $E_i$. Finally, 
$\sigma_2^{(i)} = ((\sigma_0)^{(i)})^{-1} \sigma_1^{(i)}$.
The genus defect function is the restriction of the
genus defect function of $G$ to the subgraph. 

Now we can
introduce a structure of a stable ribbon graph on the set
of remaining edges $E_r$. The underlying graph is obtained
from $G$ by contracting the edges of $E_c$. The permutation
$\sigma_1'$ is the restriction of $\sigma_1$ to the half-edges
of $E_r$. The image of a half-edge $h$ under $\sigma_2'$
is the first half-edge among $\sigma_2(h)$, $\sigma_2^2(h)$, \dots
to belong to an edge of $E_r$. Finally, 
$\sigma_0'= \sigma_1' \sigma_2'^{-1}$. If a vertex of the
new graph is the result of the contraction of $E_i$, then
its genus defect is equal to the genus of the stable
ribbon graph $E_i$. The genus defects of the other vertices
are the same as they were in the graph $G$.

\begin{proposition} \label{Prop:contraction}
The stable ribbon graph $E_r$ described above is the
result of a contraction of the edges $E_c$ in the 
stable ribbon graph $G$.
\end{proposition}

Now, using stable ribbon graphs, we can
construct a new orbifold cell complex.

\begin{definition} \label{Def:A}
Denote by $A$ the following orbifold cell complex. Its cells
are in a one-to-one correspondence with stable
ribbon graphs with $n$ numbered faces. 
If $G$ is such a graph and $E$ the set of its edges,
the corresponding cell $\Delta$ is isomorphic to $\R_+^E/\Aut(G)$,
where $\Aut(G)$ is the (finite) group of automorphisms
of $G$. A cell $\Delta_1$ is a boundary cell of the
cell $\Delta_2$ if the corresponding stable ribbon
graph $G_1$ can be obtained from the stable ribbon
graph $G_2$ by contracting several edges $e_1, \dots, e_k$. 
The cell $\Delta_1$ is then glued to $\Delta_2$
along $e_1 = \dots = e_k =0$.
\end{definition}

\begin{definition}
Denote the $\B_i$ the orbifold cell complex of pairs $(G,x)$, where
$G$ is a stable graph with $n$ numbered faces and $x$ a point
lying on its $i$th face $B_i$. 
\end{definition}

The fact that $\B_i$ indeed has a natural structure
of a cell complex can be seen in the following way.

Consider a cell $\Delta$ of $A$. Its preimage in
the total space of $\B_i$ (under the projection to $A$
forgetting the point $x$) 
is naturally subdivided into cells. Each cell is
composed of the points $(G,x)$ such that $x$ 
lies on some given edge of the polygon $B_i$, or such that
$x$ coincides with one of the vertices of $B_i$.
These cells are then glued to each other in the
obvious way.

\subsection{Factorizing $\MM_{g,n}$}

In Section~\ref{Ssec:stabStreb} we saw that a Strebel
differential on a stable curve vanished identically
on the irreducible components that do not contain
marked points. Therefore it is a good idea to contract
such components.

Let $C$ be a stable curve with $n$ marked points. 
Consider the curve $\widetilde C$ obtained from $C$
by contracting to a point each unmarked component
(i.e., an irreducible component
of $C$ that does not contain marked points). 

On the curve $\widetilde C$ we can define a genus
defect function. It is a function with a finite support
and with positive integer values, defined in the following way.

Consider the subcurve of $C$ composed of its
unmarked components.
Each connected component of this subcurve is
contracted to a point of $\widetilde C$ and we assign
to this point the arithmetic genus of the corresponding
contracted component (cf. Definition~\ref{Def:stabgraph}).

\begin{definition} \label{Def:contraction}
We call the curve $\widetilde C$ endowed with the
genus defect function the {\em contraction} of $C$.
\end{definition}

\begin{definition}\label{Def:MR}
We call the {\em minimal reasonable} compactification
of $\M_{g,n}$ (denoted by $K\MM_{g,n}$)
the quotient of the Deligne-Mumford compactification 
$\MM_{g,n}$ by the following equivalence relation: two points 
$x,y \in \MM_{g,n}$ are equivalent if the corresponding 
contractions ${\widetilde C}_x$ and ${\widetilde C}_y$ 
are isomorphic.
\end{definition}

This compactification was defined in Kontsevich's original
paper~\cite{Kontsevich}. Looijenga (\cite{Looijenga95}, Lemma~3.1)
shows that it is a compact Hausdorff topological orbifold.
It is not known whether it can be given a natural algebraic
structure.

Note that we have constructed $K\MM_{g,n}$ 
together with a projection
$$
\MM_{g,n} \rightarrow K\MM_{g,n}.
$$
This projection contracts some subvarieties of $\MM_{g,n}$
of complex codimension at least $1$. Therefore the fundamental
homology class of $\MM_{g,n}$ is sent to the fundamental
homology class of $K\MM_{g,n}$.

\begin{proposition}
The line bundles $\L_i$ over $\MM_{g,n}$ are pull-backs of
some complex line bundles over $K\MM_{g,n}$, that we will also denote by
$\L_i$.
\end{proposition}

\paragraph{Proof.} This is almost obvious. Indeed, if two
curves $C_1$ and $C_2$ are equivalent in the sense of
Definition~\ref{Def:MR}, then the cotangent lines $L_i$
to marked points on both curves are naturally identified.
\qed

\bigskip

\begin{proposition}
The intersection numbers of the first Chern classes
$c_1(\L_i)$ are the same on any compactification $X$
of $\M_{g,n}$ that can be projected on $K\MM_{g,n}$ in such a way
that the line bundles $\L_i$ are obtained by pull-back from
$K\MM_{g,n}$ and the fundamental class of $X$
is sent to the fundamental class of $K\MM_{g,n}$.
\end{proposition}

\paragraph{Proof.} This is again obvious.
Instead of calculating the intersection
numbers on $X$ we can calculate them on $K\MM_{g,n}$ 
and pull them back on $X$. \qed

In particular all the intersection numbers are the
same on $\MM_{g,n}$ and $K\MM_{g,n}$

\bigskip

It is not known whether
the compactifications $K\MM_{g,n}$ can be endowed
with the structure of singular algebraic varieties.
It had been conjectured that this can be achieved using
the semi-ampleness of the line bundles $\L_i$ over
$\MM_{g,n}$. However, Sean Keel showed in~\cite{Keel},
Section~3 that the semi-ampleness fails for $g \geq 3$
(although in finite characteristic it holds for any $g$ and $n$,
see~\cite{Keel2}).

For $g=0$, the algebraic version of the minimal reasonable
compactification $K\MM_{0,n}$ was constructed by Marco
Boggi~\cite{Boggi}. He discribed it as a solution to a
moduli problem using the construction that we give
below. He also gave a description of $K\MM_{0,n}$
using blow-ups of a projective space and studied the
action of the symmetric group on it.

The space $K\MM_{0,n}$ was also defined and used
in~\cite{GorLan} V.~Goryunov and S.~Lando, although
the authors did not know its exact interpretation as 
a moduli space.

Let us sum up the algebraic construction of $K\MM_{0,n}$.
Consider a smooth rational curve with $n$ marked points
$(\CP^1, x_1, \dots, x_n) \in \M_{0,n}$.
Consider a degree~$1$ rational function $f_i$ on $\CP^1$ 
having its pole at $x_i$ and such that
$$
f_i(x_1) + \dots + f_i(x_{i-1}) + 
f_i(x_{i+1}) + \dots + f_i(x_n) = 0.
$$
Such a function is unique, up to a multiplicative
constant, therefore its values 
$$
\biggl( f_i(x_1),  \dots, f_i(x_{i-1}),
f_i(x_{i+1}), \dots, f_i(x_n) \biggr)
$$
yield a map $F_i : \M_{0,n} \rightarrow \CP^{n-3}$.
Putting together such maps for all $i$,
we obtain a map $F : \M_{0,n} \rightarrow (\CP^{n-3})^n$.
We claim that the closure of the image of $F$
in $(\CP^{n-3})^n$ is an algebraic model of the
minimal reasonable compactification. More precisely,
$F$ can be naturally extended to a holomorphic map 
$F: \MM_{0,n} \rightarrow (\CP^{n-3})^n$ that sends
two stable curves to the same point if and only if
their contractions (Definition~\ref{Def:contraction})
are isomorphic.

Let us sketch the proof of this fact. First of all, the
family of functions $f_i$ can be extended to a family of
functions on the compactification $\MM_{0,n}$. The function
$f_i$ on a stable curve $C$ with $n$ marked points is defined
as follows. It is constant on each irreducible component
of $C$ that does not contain $x_i$. On the component that
contains $x_i$, the function $f_i$ is of degree~$1$ and
has a simple pole at $x_i$. And, as before, we have
$$
f_i(x_1) + \dots + f_i(x_{i-1}) + 
f_i(x_{i+1}) + \dots + f_i(x_n) = 0.
$$
Such a function $f_i$ is, again, unique, up to a multiplicative
constant. It is easy to see that the values of the function
$f_i$ at the points $x_j, j \not= i$ allow one to reconstitute
the complex structure on the component of $C$ that contains
$x_i$, but not on the other components.

Thus the function $F$ can be extended to $\MM_{0,n}$ and a
point $F(C)$ allows one to reconstitute the complex structure
of the marked components of $C$, but not that of the unmarked
components.

\subsection{A homeomorphism between $K\MM_{g,n} \times \R_+^n$
and the cell complex of stable ribbon graphs}

Consider two stable curves $C_1$ and $C_2$ that
are mapped to the same point of $K\MM_{g,n}$. Let
$p_1, \dots, p_n$ be a given list of perimeters (positive
real numbers). It is clear that the stable ribbon graphs
formed by the nonclosed horizontal trajectories of the
Strebel differentials on $C_1$ and $C_2$ are the same,
including the edge lengths. We can therefore define a map 
$h$ from $K\MM_{g,n} \times \R_+^n$
to the cell complex $A$ of stable ribbon graphs with
$n$ numbered faces and endowed with edge lengths
(see Definition~\ref{Def:A}).

\begin{theorem}\label{Thm:homeomorphism}
The map
$$
h: K\MM_{g,n} \times \R_+^n \rightarrow A
$$
is an isomorphism of topological orbifolds. The
polygonal bundles $\B_i$ over $A$ are naturally
identified with the real spherizations 
$(\L_i^* \setminus \mbox{\rm zero section})/\R_+$ of the complex line bundles
$\L_i^*$ dual to $\L_i$.
\end{theorem}

This theorem was formulated by M.~Kontsevich without a
proof. It follows from the main theorem (Theorem~8.6)
of E.~Looijenga's paper~\cite{Looijenga95}. Here we give
a different proof.

\paragraph{Proof of Theorem~\ref{Thm:homeomorphism}.}
The identification of $(\L_i^* \setminus \mbox{\rm zero section})/\R_+$
with $\B_i$ is immediate (see the discussion in the end
of Section~\ref{Ssec:stabStreb}).

The bijectivity of $h$ is a reformulation of
Proposition~\ref{Prop:strips2}.

The continuity of $h$ and that of $h^{-1}$ are
equivalent. Indeed,
both spaces $K\MM_{g,n} \times \R_+^n$ and 
$A$ have natural proper projections to $\R_+^n$.
We know that the map $h$ is a bijection that commutes
with the projections. Therefore if $h$ or its inverse
is continuous, then $h$ is a homeomorphism. 

Thus the main task is to prove the continuity of $h$, which
we will do using Theorem~\ref{Thm:cont2} and 
Proposition~\ref{Prop:teichmetric}.

{\bf 1. A sequence of Strebel differentials.}

Consider a sequence of stable curves $C_m$ tending to
a stable curve $C$, together with a sequence of
$n$-tuples $p^{(m)} = (p_1^{(m)}, \dots, p_n^{(m)})$
of positive real numbers tending to an $n$-tuple
$p = (p_1, \dots, p_n)$ of positive real numbers.

From Theorem~\ref{Thm:cont2} we know that in
the vector bundle $\ocW$ of quadratic differentials,
the sequence of the corresponding Strebel differentials $\varphi_m$
on $C_m$ tends to the Strebel differential $\varphi$
on $C$. Moreover, according to Proposition~\ref{Prop:teichmetric},
there is a sequence of deformations $f_m:C_m \rightarrow C$
such that the sequence $(f_m)_* \varphi_m$
converges to $\varphi$, uniformly on any compact set 
$K \subset (C \setminus \mbox{marked poles and nodes})$. 

We must prove that in the cell complex $A$,
the sequence of stable graphs with edge lengths corresponding to
$\varphi_m$ tends to the stable graph with edge lengths
corresponding to $\varphi$.

{\bf 2. For $m$ big enough, a disc domain of $\varphi_m$ contains
any compact set inside the corresponding disc domain of $\varphi$.}

Let $z=z_i$ be one of the marked points in the curve $C$,
let $D$ be the corresponding disc domain consisting of closed
horizontal trajectories, and let $K \subset D$ be any compact
set. Then, for $m$ big enough, $f_m^{-1}(K)$ is contained
in the disc domain of $\varphi_m$ surrounding the marked
point $f_m^{-1}(z)$. Indeed, consider a compact 
annulus $K'$ surrounding $K$. We suppose that $K'$ is 
composed of closed trajectories of the disc domain $D$
and that $K \cap K'= \emptyset$. The annulus $K'$ does
not contain nodes of $C$, nor zeroes or poles of $\varphi$.

In the proof of Theorem~\ref{Thm:contfam}, paragraph {\bf 2},
we proved that, for $m$ big enough, 
$f_m^{-1}(K')$ necessarily contains a closed horizontal
trajectory of $\varphi_m$ that makes exactly one turn
around $f_m^{-1}(K')$. This trajectory surrounds $f_m^{1}(K)$
entirely. Thus $f_m^{-1}(K)$ lies inside a disc
domain of $\varphi_m$.

{\bf 3. Cutting the curve $C$ into pieces.}

For shortness, we call a {\em disc neighborhood} of
a point $x \in C$ any neighborhood of $x$
homeomorphic to an open disc.

To every marked point $z_i$ on $C$ we assign a
disc neighborhood $U_i \subset C$.

To every vertex $v$ of the stable graph of the differential
$\varphi$ we assign an open set $U_v \subset C$ as follows.
If $v$ is a zero of $\varphi$ at a regular point of
$C$, then $U_v$ is just a disc neighborhood of $v$.
If $v$ is a node of $C$ such
that both irreducible components intersecting at $v$ are marked
(contain marked points), then $U_v$ is the union of
two disc neighborhoods of $v$ on both components.
If $v$ is obtained by contracting one or more unmarked
irreducible components of $C$, then $U_v$ is the union
of these unmarked components and of the disc neighborhoods
of all the nodes at which they meet marked components.

Finally, to every edge $e$ we assign a compact set $K_e \subset C$.
It is homeomorphic to a closed disc and contains in its
interior the part of $e$
that lies outside $U_v$ and $U_{v'}$, where $v$ and $v'$
are the vertices of $e$. Moreover, we suppose that $K_e$
does not intersect the vertical trajectories of $\varphi$
that join marked points to the vertices $v$ and $v'$.

\begin{figure}
\begin{center}
\
\epsfbox{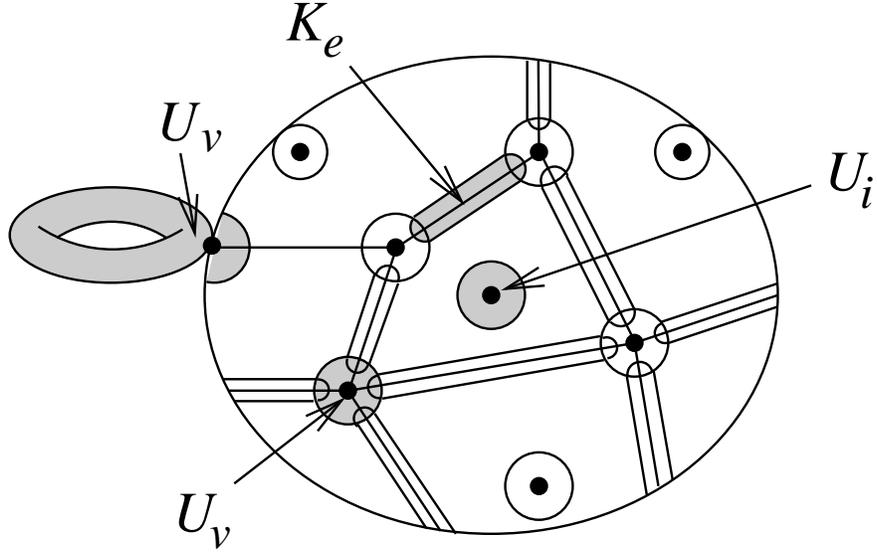}

\caption{The sets $U_i$, $U_v$, and $K_e$.}
\label{Fig:pieces}

\end{center}
\end{figure}

All these sets are represented in Figure~\ref{Fig:pieces}.

All the disc neighborhoods above can be chosen
arbitrarily small.

We denote by $G = \Gr(\varphi)$ the stable ribbon graph endowed
with edge lengths, assigned to the Strebel differential $\varphi$.
Similarly, $G_m = \Gr(\varphi_m)$ is the stable ribbon
graph with edge lengths assigned to $\varphi_m$.

{\bf 5. For $m$ big enough, the vertices of $G_m$ in $C_m$
lie inside the union of the sets $f_m^{-1}(U_v)$ over all
the vertices $v$ of $G$.}

Indeed, consider the compact set $K \subset C$ obtained
from $C$ by taking away all the open sets $U_i$ and $U_v$.
There are no zeroes or poles of $\varphi$ nor nodes of $C$
in $K$. Therefore, for $m$ big enough, there are no
vertices of $G_m$ in $f_m^{-1}(K)$. On the other hand,
according to {\bf 3}, there are no vertices of $G_m$ in
$f_m^{-1}(U_i)$, because $f_m^{-1}(U_i)$ belongs
to the $i$th disc domain of $\varphi_m$. Thus every 
vertex of $G_m$ lies inside $f_m^{-1}(U_v)$ for some
vertex $v$ of $G$.

Conversely, let us prove, that each set $f_m^{-1}(U_v)$
contains at least one vertex of $G_m$. 

Suppose $v$ is a $k$-tuple zero of $\varphi$ at a regular point 
of $C$. Then we consider a circle surrounding $v$ and lying
inside $U_v$. As we go around this circle, the horizontal
trajectories of $\varphi$ make $-k$ half-turns with respect to the
tangent line to the circle. Therefore the same is true for
the horizontal trajectories of $\varphi_m$, for $m$ big enough.
Thus $f_m^{-1}(U_v)$ contains one or several zeroes of $\varphi_m$
whose sum of multiplicities equals $k$.

Suppose $U_v$ contains at least one node of $C$. If $f_m^{-1}(U_v)$
also contains a node of $C_m$ there is nothing to prove, because
this node is necessarily a vertex of $G_m$. Suppose that
$f_m^{-1}(U_v)$ does not contain nodes. Then it is a
smooth Riemann surface with holes. We must prove that
$\varphi_m$ has at least one zero on this surface.

First of all, if we are given a quadratic differential
without poles on a Riemann surface $S$ of genus
$g$ with $d$ holes,
the number of its zeroes (taking their multiplicities
into account) equals
$$
|\mbox{zeroes}| = 4g - 4 + 2d - 2 |\mbox{turns}|,
$$
where $|\mbox{turns}|$ is the total number of turns that the 
horizontal trajectories make with respect to the tangent
lines to the boundary circles.

Recall that $U_v$ contains several unmarked components
of $C$ and several small discs surrounding nodes
of marked components. 
As above, we draw a circle inside each such small disc,
surrounding the corresponding node. If the node is a
vertex of degree $k$, then, as
we go around the circle, the
horizontal trajectories make $-k$ half-turns
with respect to the tangent lines to the circle. 
Therefore the same is true for $(f_m)_*\varphi_m$,
for $m$ big enough. Since
$f_m^{-1}(U_v)$ has either a positive genus
or at least two holes, and since the
number of turns is always negative, $\varphi_m$ has
at least one zero in $f_m^{-1}(U_v)$.

{\bf 6. An injection from the set of edges of $G$ to the
set of edges of $G_m$.}

Let $e$ be an edge of $G$. We will assign to $e$ an edge
of $G_m$. Later we will see that one obtains $G$ by 
contracting the edges of $G_m$ that are not assigned to any edge 
of $G$.

A vertical trajectory of a Strebel differential usually joins
two marked points (that can happen to be the same) and
crosses exactly one nonclosed horizontal trajectory.
(The exceptions are those vertical trajectories
that join a marked point to a zero of the differential
or to a node of the curve.)

Two vertical trajectories intersect the same nonclosed
horizontal trajectory if and only if they join the same
pair of marked points and, moreover, bound a region in the
stable curve, that does not contain zeroes of poles
of the Strebel differential, nor nodes of the curve.

Consider the vertical trajectory $\alpha$ of $\varphi$
through any point $x \in K_e$. By the choice of $K_e$,
we know that $\alpha$ does not end at a vertex of $G$.
Therefore it joins some marked points $z_i$ and $z_j$
(it can happen that $z_i = z_j$). Denote
by $\alpha_m$ the vertical trajectory of $\varphi_m$
through $f_m^{-1}(x)$. For $m$ big enough,
$f_m(\alpha_m)$ follows $\alpha$ closely enough to enter the
neighborhoods $U_i$ and $U_j$. Therefore it necessarily
joins $z_i$ and $z_j$, just as $\alpha$ does (because
$f_m^{-1}(U_i)$ and $f_m^{-1}(U_j)$ lie inside the
corresponding disc domains). According to the above
remark, the vertical trajectory $\alpha_m$ crosses exactly
one nonclosed horizontal trajectory of $\varphi_m$,
in other words, exactly one edge of $G_m$. We
denote this edge by $e_m$ and assign it to the edge $e$ of $G$. 

For $m$ big enough, the resulting edge does
not depend on the choice of $x \in K_e$.
Indeed, two vertical trajectories of $\varphi_m$ 
through two points of $f_m^{-1}(K_e)$ 
bound a region in $C_m$, that does not contain
zeroes or poles of $\varphi_m$ or nodes of $C_m$. 
Therefore these two vertical trajectories cross
the same edge of $G_m$.

Let us prove that $e \mapsto e_m$ is an injection.
Consider two vertical trajectories of
$\alpha_m$ through $x \in f_m^{-1}(K_e)$ and $\alpha_m'$ through
$x' \in f_m^{-1}(K_{e'})$, where $e$ and $e'$ are two different
edges. Let us prove that they cannot intersect the same
edge of $G_m$. If they cross the same edge of $G_m$,
it means that they join the same
pair of marked points and, moreover, 
bound a region in $C_m$, that does
not contain zeroes or poles of $\varphi_m$ or nodes of $C_m$.
It is easy to see that, if $\alpha_m$ and $\alpha_m'$ join
the same pair of marked points and bound a region in $C_m$, 
then this region contains at least one set $f_m^{-1}(U_v)$ for some
vertex $v$ of $G$. But, according to {\bf 4}, this set contains
at least one vertex of $G_m$. Thus $\alpha_m$ and $\alpha_m'$
correspond to different edges of $G_m$. This proves the
injectivity of the map $e \mapsto e_m$.

Each edge $e_m$ lies entirely inside 
$f_m^{-1}(K_e \cup U_v \cup U_{v'})$, where $v$ and $v'$ are
the vertices of $e$. Indeed, $e_m$ is
entirely contained inside the union of sets $f_m^{-1}(K_*)$
and $f_m^{-1}(U_*)$ over all edges and vertices
of $G$, because the complement of this union is covered
by disc domains of $\varphi_m$. On the other hand,
$e_m$ does not meet 
$f_m^{-1}(K_{e'})$ for any edge $e' \not= e$, because otherwise
it would cross a vertical trajectory of $\varphi_m$
that is should not cross. Thus $e_m$ is contained
in $f_m^{-1}(K_e \cup U_v \cup U_{v'})$, because
only $f_m^{-1}(U_v)$ and $f_m^{-1}(U_{v'})$ 
have common points with $f_m^{-1}(K_e)$.

For the same reason, an edge of $G_m$
that does not correspond to an edge of $G$ lies
entirely in $f_m^{-1}(U_v)$ for some vertex $v$.

{\bf 7. The difference of lengths between an edge $e$ of $G$
and the corresponding edge $e_m$ of $G_m$ is less than
$\varepsilon$. The other edges of $G_m$ are shorter than
$\varepsilon$.}

At present, we have proved the following. To each edge $e$
of $G$ joining to vertices $v$ and $v'$ we can assign
an edge $e_m$ of $G_m$ joining some points inside
$f_m^{-1}(U_v)$ and $f_m^{-1}(U_{v'})$. The image
of $e_m$ under $f_m$ is contained inside the union of $K_e$,
$U_v$, and $U_{v'}$. An edge of $G_m$ that does
not correspond to an edge of $G$ lies inside
$f_m^{-1}(U_v)$ for some vertex $v$. All this, of course,
is only true starting from some $m$.

Consider an edge $e$ and let $l_e$ be its length.
If we choose the disc neighborhoods $U_v$ of the vertices small
enough, the length of the part of $e$ that lies outside
the neighborhoods $U_v$ is greater than $l_e - \varepsilon$.
Consider the line $f_m(e_m)$, more precisely, its
part lying in $K_e$. By choosing $K_e$ small enough
and $m$ big enough, we see that the the length of
$f_m(e_m) \cap K_e$ measured with the differential
$(f_m)_* \varphi_m$ differs from the length of
$e \cap K_e$ measured with $\varphi$ by less than $\varepsilon$.
Thus the length of $e_m$ is greater than $l_e - 2 \varepsilon$.

But the total sum of lengths of the edges of $G_m$ is
fixed: it is equal to the sum of the perimeters
$\sum_i p_i^{(m)}$, which is arbitrarily close
to $\sum_i p_i$. Therefore, for $m$ big enough, the length
of each edge $e_m$ is arbitrarily close to that of $e$, while
the lengths of the edges of $G_m$ that do not correspond
to edges of $G$ are arbitrarily small.

{\bf 8. The genus defect function.}

The genus defect assigned to a vertex $v$ of $G$
is equal to the arithmetic genus of the open set $U_v$, which
is actually a singular noncompact complex curve. Using
Proposition~\ref{Prop:contraction}, it is easy to see
that this genus defect is indeed obtained by
contracting the edges of $G_m$ that lie in $f_m^{-1}(U_v)$.

{\bf 9. Conclusion.}

Thus, for $m$ big enough, the stable ribbon graph $G$ is
obtained from $G_m$ by contracting some edges of length
less than $\varepsilon$ and by changing the lengths of
the other edges by less than $\varepsilon$. This means
that the graph $G_m$ lies in an $\varepsilon$-neighborhood
of $G$ in the cell complex $A$. Since $\varepsilon$ can
be chosen arbitrarily small, the sequence $G_m$ tends
to $G$ in $A$.
\qed

\subsection{Looijenga's results}

This section is a very brief review of Looijenga's
paper~\cite{Looijenga95}. We follow, as closely as possible, 
the notation introduced there.

Looijenga's main result is the continuity of a map similar
to the map $h^{-1}$ in our Theorem~\ref{Thm:homeomorphism}.
In other words, he proves that when one changes continuously 
the stable ribbon graph with edge lengths, the corresponding
Riemann surface glued from strips also changes continuously.
The main problem is that if we consider a sequence of ordinary
ribbon graphs $G_m$ converging in $A$, the corresponding
sequence of smooth curves does not necessarily converge
in the Deligne-Mumford compactification $\MM_{g,n}$, but
only in the quotient $K\MM_{g,n}$. However, there is no simple
criterion of convergence in $K\MM_{g,n}$.
To solve this problem, Looijenga
constructs a more complicated cell complex, that turns out
to be homeomorphic to $\MM_{g,n} \times \R_+^n$. This is done
roughly as follows. Consider a stable ribbon graph $G$
and normalize its edge lengths (by multiplying
them by a constant) so that their sum
equals $1$. Choose a subset $E_1$ of
the set of edges $E$, and suppose the lengths of the edges
of $E_1$ tend to $0$. Instead of forgetting everything
about the lengths of the edges of $E_1$ (as we do when
we contract them to obtain a new stable ribbon graph), we
normalize their lengths anew, so that their sum equals $1$.
Among the edges of $E_1$ their can now be a new subset
of edges $E_2$ whose lengths still tend to $0$. We 
normalize them once again, so that their sum equals $1$.
And so on. Using all this information, we can construct
a cell complex with a projection onto $A$, and such that
a converging sequence in this complex induces a converging
sequence of stable curves.

Looijenga's construction is a little more general than what is
needed for Kontsevich's proof, because he allows the
perimeters $p_1, \dots, p_n$ to vanish (but under the above
normalization their sum remains equal to $2$, so at least
one perimeter must remain positive) and studies what
happens to Strebel differentials and
stable ribbon graphs in that case. Therefore his definitions
of stable graphs and of the minimal reasonable compactification
are slightly different from ours. 

Suppose the set of marked
points $\{ 1, \dots, n \}$ is divided into two disjoint
parts: $V \sqcup Q$, $Q \not= \emptyset$. Here $V$ is
the set of points such that the corresponding perimeters vanish,
while $Q$ is the set of points such that the corresponding
perimeters do not vanish.

\paragraph{Stable ribbon graphs.} The natural modification
of the notion of Strebel differentials to this case is to
consider the Strebel differentials on our surface
punctured at the points of the set $V$ (and with double
poles with given residues at the points of $Q$). It is easy
to see that if one puts the points of $V$ back into the
surface, the Strebel differential will have at most simple
poles at these points. Therefore they will be vertices of
the graph of nonclosed horizontal trajectories. Thus
the new stable graphs, instead of having $n$ numbered
faces, have $n$ numbered faces or vertices. We do not
give the precise definition of a stable ribbon graph
in this setting; it is a formalization of the properties
of the graph of nonclosed horizontal trajectories
of a Strebel differential.

\paragraph{Minimal reasonable compactifications.} In our
definition of the minimal reasonable
compactification, two points $x,y \in \MM_{g,n}$
are identified if the stable curves $C_x$ and
$C_y$ give the same curve when one contracts
each of their components that do not contain marked points.
Analogously, $K_Q\MM_{g,n}$ is the quotient of
$\MM_{g,n}$ in which two points $x,y \in \MM_{g,n}$
are identified if the stable curves $C_x$ and
$C_y$ give the same curve when one contracts
each of their components that do not contain 
the marked points {\em of the set $Q$}.

\paragraph{Teichm\"uller spaces.} Looijenga works with
Teichm\"uller spaces rather than moduli spaces. The
advantage of this approach is that the spaces considered
are not orbifolds, but usual topological space. However
one has to work with non locally compact topological spaces.

Let ${\cal T}_{g,n}$ be the Teichm\"uller
space of Riemann surfaces of genus $g$ with $n$ marked points. Its
quotient by the action of the mapping class group $\Gamma$
is the moduli space $\M_{g,n}$. \,
Looijenga uses an augmented Teichm\"uller space $\T_{g,n}$ constructed
by Harvey~\cite{Harvey}. The space $\T_{g,n}$ is endowed with a proper 
action of the mapping class group $\Gamma$ and there is
a natural $\Gamma$-invariant surjective projection from
$\T_{g,n}$ onto the Deligne-Mumford compactification $\MM_{g,n}$. 

The space $\T_{g,n}$ is constructed as follows.
Let $C$ be a smooth genus $g$ Riemann surface with $n$ punctures. 
We can choose $3g-3+n$ simple loops in $C$ in such a way that they
cut $C$ into $2g-2+n$ ``pants'' ($3$-holed spheres). 
In the free homotopy
class of each loop we can choose the shortest geodesic
with respect to the unique complete metric of curvature $-1$
on $C$, compatible with the conformal structure. 
To each geodesic we can assign its length $l \in \R_+$
and the angle $\theta \in \R$,
with respect to some chosen gluing, at which the two 
pants adjacent to it are glued to each other.
These lengths and angles are called the {\em Fenchel-Nielsen
coordinates} in the Teichm\"uller space.
It is well-known that they determine a real analytic 
diffeomorphism of the Teichm\"uller space onto the open
octant $\R_+^{3g-3+n} \times \R^{3g-3+n}$
(see, for example,~\cite{Abikoff}). Now we simply add
to the octant the boundary hyperplanes and endow the obtained
closed octant with the usual topology. This corresponds to pinching
the geodesics, but retaining the angles at which the
adjacent pants were glued to each other. This operation
can be carried out for all possible choices of 
$3g-3+n$ geodesics and the points we adjoin to the Teichm\"uller
space add up into the augmented Teichm\"uller space $\T_{g,n}$.
For a smooth surface $S$ with $n$ marked points we
define a {\em stable complex structure} to be a complex
structure  $J$ defined on $S \setminus L$, where
$L \subset (S \setminus \mbox{marked points})$ is a finite
set of simple loops, such that pinching the loops we
obtain a stable curve. Then, 
as a set, $\T_{g,n}$ can be identified with the set of all stable
complex structures on a given smooth surface $S$ with
$n$ marked points, up to diffeomorphisms homotopic
to the identity, relatively to the marked points.

Looijenga uses the following criterion of convergence
on $\T_{g,n}$. Let $J_m$ be a sequence of complex structures on 
a smooth surface $S$ with $n$ marked points
and let $J$ be a stable complex structure on $S \setminus L$.
If $J_m$ converges to $J$ uniformly on every compact set 
$K \subset (S \setminus L)$, then the sequence of curves
$(S, J_m)$ converges to $(S, J)$ in the space $\T_{g,n}$.

Starting from $\T_{g,n}$, one can construct the analog of minimal reasonable
compactifications for Teichm\"uller spaces. Indeed,
there is a $\Gamma$-invariant map $\T_{g,n} \rightarrow K_Q\MM_{g,n}$
for any subset $Q$ of the set of marked points.
We denote by $K_Q {\cal T}_{g,n}$ the space obtained from $\T_{g,n}$ 
by contracting to one point each connected component 
of the preimage of each point under this map. This space
is still endowed with an action of $\Gamma$, although it is
not proper any longer.

\paragraph{The simplicial complex of stable graphs.}
Let $S$ be a fixed surface with $n$ marked points.
Consider the following infinite simplicial complex $A_{\cal T}$.

Consider the set of homotopy classes (relatively to
the marked points) of simple non-oriented
arcs in $S$ joining two (possibly coinciding)
marked points and avoiding the other marked points.
The homotopy class is called trivial if the
corresponding arc is contractible in $S$ punctured
at the marked points.

A vertex of $A_{\cal T}$ is a nontrivial homotopy class as
above. A set of homotopy classes forms a simplex if 
the homotopy classes can by realized by loops that
do not intersect (except at their endpoints).

There is a natural action of the mapping class group $\Gamma$
on $A_{\cal T}$. The quotient of $A_{\cal T}$
by this action is a finite orbifold simplicial complex.
One can prove that the quotient complex $A_{\cal T}/\Gamma$
is isomorphic to the complex $A$ defined by stable ribbon graphs. 
(If a stable graph is a just a ribbon graph with $n$ numbered
faces and with degrees of vertices $\geq 3$, then the
corresponding set of arcs is obtained by considering
the dual graph: joining by arcs the centers of adjacent faces.
When we contract an edge in the stable ribbon graph
we must erase the corresponding arc.)

Thus there are two equivalent
ways of defining the same simplicial complex $A$. The definition
with isotopy classes of arcs is certainly more elegant,
but stable ribbon graphs are needed anyway to make
the connection with Strebel differentials.

\paragraph{A quotient of \,  $\T_{g,n} \times \mbox{ simplex}$.}

Finally, let $\Delta_n$ be the standard $n$-simplex.
We consider the topological space $|K_{\bullet} {\cal T}|$
obtained from $\T_{g,n} \times \Delta_n$ by the following
factorization. Consider a point $x$ in $\Delta_n$ and let
$Q$ be the set of its nonzero coordinates
(a subset of $\{ 1, \dots, n\}$). Then the ``layer''
$\T_{g,n} \times \{ x \}$ is factorized
so as to obtain $K_Q {\cal T} \times \{ x\}$.

\begin{theorem} {\rm \cite{Looijenga95}} \label{Thm:looijenga}
There is a natural bijective continuous map from the complex
$A_{\cal T}$ to $|K_{\bullet} {\cal T}|$. It is $\Gamma$-invariant
and commutes with the projections of both spaces
on the simplex $\Delta_n$.
\end{theorem}

The above mapping is not a homeomorphism, but if we quotient
both spaces by the action of $\Gamma$ it becomes a 
homeomorphism, because a continuous bijection between
two compact topological spaces is necessarily a homeomorphism.
If, in both spaces, we take the preimage of the interior
of the simplex $\Delta_n$ we immediately obtain 
Theorem~\ref{Thm:homeomorphism}.

Looijenga proves Theorem~\ref{Thm:looijenga} by constructing
an explicit trivialization of families of Riemann surfaces
as a ribbon graph tends to a stable ribbon graph and
using the convergence criterion that we formulated in the
paragraph on Teichm\"uller spaces.

\section{The Chern classes $c_1(\L_i)$} \label{Sec:Chern}

Here we recall Kontsevich's expression for the
first Chern classes $c_1(\L_i)$. Theorem~\ref{Thm:homeomorphism}
allows us to work on the cell complex $A$.\, Kontsevich's
expressions are cellwise smooth continuous differential
forms, and we explain the framework in which such forms
can be used.

\subsection{A connection on the bundles $\B_i$}
\label{Ssec:conncurv}

Once we have found a homeomorphism $h$ between
$K\MM_{g,n} \times \R_+^n$ and the cell complex $A$,
we will no longer use the smooth
structure of $K\MM_{g,n}$ (defined outside the
singularities). Instead, we use the natural piecewise
smooth structure of the cell complex $A$. The relation
between the two is rather delicate and we won't discuss
it here.

Thus we have $n$ polygonal bundles $\B_i$ over the cell
complex $A$ and we want to find their first Chern
classes.

(We remind the reader how to define the first Chern
class of a topological oriented circle bundle over any topological 
space $X$ homotopically equivalent to a cell complex.
There exists a continuous map $f$ from $X$ to the infinite
projective space $\CP^{\infty}$ such that the circle
bundle over $X$ is isomorphic to the pull-back under $f$ of
the canonical circle bundle over $\CP^{\infty}$. The
first Chern class of the bundle is the pull-back
under $f$ of the natural $2$-cohomology class of 
$\CP^{\infty}$. It is an element of $H^2(X, \Z)/\mbox{torsion}$.)

Consider one of the polygonal bundles $\B= \B_i$. 
Kontsevich constructs an explicit $1$-form $\alpha$ on each
cell of the total space of $\B$, claiming
that $d \alpha$ represents the first Chern class
of the line bundle $\L_i$ over $\MM_{g,n} \times \R_+^n$
(see~\cite{Kontsevich}, Lemma~2.1).
The $1$-form in question is the following.

Let $p$ be the perimeter of the polygon $B$, $k$ its
number of vertices, and 
$$
0 \leq \phi_1 < \dots < \phi_k < p
$$
the distances from the distinguished point of the
polygon to its vertices
(as we go around the polygon counterclockwise). 
Moreover, denote by $l_i$, $1 \leq i \leq k$, the length of the
edge that follows the $i$th vertex. Then we have
$$
\alpha = \sum_{i=1}^k \frac{l_i}p \; d \! 
\left( \frac{\phi_i}p \right) \, .
$$

\subsection{Differential geometry on polytopal complexes}

Here we introduce a framework for working with 
cellwise smooth differential forms on cell complexes.
It has appeared, for example, in Sullivan's work~\cite{Sullivan}, 
but we give a complete exposition here, adapted to our needs.

The results of this section can be considered as
a far-reaching generalization of the fact that
the Newton-Leibniz formula 
$\int_a^b f'(t) dt = f(b) - f(a)$ holds not only
for differentiable functions $f$, but also for
continuous piecewise differentiable ones.

First we define polytopal complexes, which are 
simply spaces glued from affine polytopes.

\begin{definition}
A {\em polytope} in a real vector space is an
intersection of a finite number of open or closed
half-spaces such that its interior is non-empty.
Replacing, in the above intersection,
some of the closed half-spaces by their
boundary hyperplanes, we obtain a {\em face}
of the polytope.
\end{definition}

Thus a polytope is always convex, but not necessarily
closed or bounded. A face is a subset of the polytope.

\begin{definition}
A {\em polytopal complex} is a finite set $X$ of polytopes
in real vector spaces, together with gluing functions
satisfying the following conditions. (i)~Each gluing
function is an affine map that identifies a polytope
$P_1 \in X$ with a face of another polytope $P_2 \in X$. 
(For brevity, we will say that $P_1$ is a face of $P_2$.)
(ii)~If $P_1$ is a face of $P_2$, which is a face of
$P_3$, then $P_1$ is a face of $P_3$, and the corresponding
gluing functions form a commutative diagram. (iii)~If $P_1 \in X$
is identified with a face of $P_2 \in X$, no other polytope
$P_1' \in X$ can be identified with the same face of $P_2$.
\end{definition}

Now we define differential forms on polytopal complexes.
A differential form on a polytope is simply a differential
form with smooth coefficients defined in some neighborhood
of the polytope in the ambient vector space.

\begin{definition} \label{Def:forms}
A {\em differential $k$-form} on a polytopal complex is
a set of differential $k$-forms defined on all the polytopes
such that restricting the $k$-form to
a face of a polytope coincides with the $k$-form on 
the face.
\end{definition}

\begin{example}
Consider two squares lying in the half-planes $x \geq 0$
and $x \leq 0$ and having a common side on the $y$ axis.
They form a polytopal complex. A differential $1$-form
on this complex can be given, for example, by
$dx+ dy$ in the right-hand square, $-2dx+dy$ in the
left-hand square, and $dy$ on their common edge. Moreover,
we could have added to our complex a third square
(or another polygon) such that the three of them
would share a common edge. The $1$-form can then
be extended to this new polygon.
\end{example}

This example shows that $k$-forms in two adjacent polytopes
sharing a common face of dimension $\geq k$ 
are not independent: they must coincide on the common face.

\begin{definition}
The exterior product and the differential $d$ of 
differential forms on polytopal complexes are defined
polytope-wise.
\end{definition}

It is obvious that we have 
$$
d^2\alpha = 0 \quad \mbox{and} \quad
d(\alpha \wedge \beta) = (d \alpha) \wedge \beta
+ (-1)^{\deg\alpha}\alpha \wedge (d \beta),
$$
because these identities are true on each polytope.
Therefore each polytopal complex possesses a de Rham
complex and the de Rham cohomology forms an algebra.

\begin{proposition}
The de Rham cohomology groups of a polytopal complex $X$ are
canonically identified, as real vector spaces, 
with its usual cohomology groups over $\R$..
\end{proposition}

\paragraph{Proof.} The de Rham complex can be
considered as a complex of sheaves on the polytopal
complex $X$. It suffices to prove that it is a
flasque resolution of the constant sheaf $\R$ on $X$.
In other words, we must prove that locally each
closed differential form is exact (except for the constant
functions considered as $0$-forms). Consider a point
$x \in X$ and a closed $k$-form $\alpha$ defined in
a sufficiently small neighborhood $U$ of $x$. We will construct
a $(k-1)$-form $\beta$ in $U$, such
that $d\beta = \alpha$. The value of $\beta$ on $k-1$ vectors
tangent to one of the polytopes is obtained by the following
standard procedure. We construct on the $k-1$ vectors
a small parallelepiped $P$ that fits entirely into the polytope.
Then we consider the cone with vertex $x$ and with base $P$.
The integral of $\beta$ over $P$ is, by definition, equal
to the integral of $\alpha$ over the cone. By letting the
sides of $P$ tend to $0$ we find the value of $\beta$
on the $k-1$ vectors. It is obvious that $\beta$ is
a $(k-1)$-form on the polytopal complex in the neighborhood
of $x$ (in the sense of Definition~\ref{Def:forms}). 
It is easy to check that if $d \alpha= 0$, then $d\beta = \alpha$.
\qed

Now the main task is to prove that the Stokes formula is
still true for differential forms on polytopal complexes.

\begin{definition}
A {\em $k$-piece} in a polytopal complex $X$ is an
affine map from a compact $k$-dimensional polytope
to a polytope of $X$.
A {\em $k$-chain} in a polytopal complex $X$ is 
a finite linear combination $C$ of $k$-pieces with
real coefficients.
The {\em boundary} of a chain and the integral of a 
$k$-form over a $k$-chain are defined in the obvious
way.
\end{definition}

\begin{proposition} {\bf (Stokes formula)}
Let $X$ be a polytopal complex, $C$ a $k$-chain
in $X$, and $\alpha$ a $(k-1)$-form on $X$. Then 
$$
\int_C d\alpha  = \int_{\d C} \alpha.
$$
\end{proposition}

\paragraph{Proof.} The formula is obvious if $C$
is composed of a unique $k$-piece, because in that case the
piece is contained in a unique polytope
of $X$. In the general case the formula is obtained
by summing over the pieces
of $C$. \qed

\begin{proposition}
The algebra structure of the de Rham cohomology of a 
polytopal complex $X$ (given by the multiplication of forms)
coincides with the usual algebra structure of the cohomology
of $X$.
\end{proposition}

\paragraph{Proof.} Recall the usual definition of
the product in the space of cohomologies (see~\cite{Masey},
chapter XIII).
Let $X$ be a polytopal complex and consider the 
polytopal complex $X \times X$ with the two projections,
$p_1$ and $p_2$, on $X$.\, For $u,v \in H^*(X, \R)$ one defines
$u \otimes v \in H^*(X \times X, \R)$ by the formula
$$
(u \otimes v) (a \times b) = (-1)^{\deg v \deg a} u(a)v(b),
$$ 
where $a$ and $b$ are two cycles in $X$. The product cycles
$a \times b$ span the whole homology group
of $X \times X$, therefore the above formula defines
the class $u \otimes v$ unambiguously. It is clear
that if differential forms $\alpha$ and $\beta$ on $X$
represent the classes $u$ and $v$, then the form 
$p_1^*\alpha \wedge p_2^*\beta$ represents the class
$u \otimes v$. Now, the product $uv$ is defined
by taking the restriction of $u \otimes v$ to the
diagonal of $X \times X$. Thus it is represented
by $\alpha \wedge \beta$. \qed

\bigskip

Now we will consider a polytopal equivalent of a
circle bundle and prove that its first Chern
classes can be expressed as the ``curvature'' of
a cellwise smooth ``connection''.

\begin{definition}
A {\em morphism} of polytopal complexes 
$F:X_1 \rightarrow X_2$ is a set $F$ of
affine maps $f:P_1 \rightarrow P_2$, where $P_1$ is
a polytope of $X_1$ and $P_2$ a polytope
of $X_2$. This set must satisfy the following
natural conditions. At least one map should
be defined on every polytope of $X_1$.
For every map $f \in F$ its restrictions to the faces
of $P_1$ must belong to $F$.
If the image of $P_1$ under $f$ belongs to a face
of $P_2$, the map from $P_1$ to this face should also
belong to $F$. If a point of $P_1$ has two different
images under maps of $F$, these images should
be identified by gluing functions of the
complex $X_2$.
\end{definition}

Note that the image of each polytope under a morphism
lies in a unique polytope of the target complex.

Let $F:X \rightarrow Y$ be a morphism of polytopal
complexes such that the preimage of each point of
$Y$ is homeomorphic to a circle. (Each such circle is
naturally subdivided into $0$-cells and $1$-cells,
and we do not require that these subdivisions be the
same for different fibers.) Suppose that the
circle bundle thus obtained is oriented. Let
$\alpha$ be a $1$-form on the polytopal complex
$X$ such that its integral over each fiber of
$F$ equals $1$ and such that $d\alpha$ is a
pull-back under $F$ of a $2$-form $\omega$ on $Y$.
The Stokes formula allows one to prove that 
$\omega$ represents the first Chern class of the
bundle. More precisely:

\begin{proposition} \label{Prop:Chern}
If $S$ is a polytopal complex homeomorphic to a
compact $2$-dimensional manifold without boundary,
and $G:S \rightarrow Y$ a morphism of complexes,
then $\int_S G^*\omega$ is equal to the first
Chern class of the pull-back to $S$ of the circle
bundle over $Y$.
\end{proposition}

\paragraph{Proof.} We can assume that $S$ is connected.
Denote by $G^*X$ the pull-back to $S$ of the bundle $X$. 
Denote by $a$ the corresponding first Chern class.
One can easily construct a section of
$G^*X$ over the surface $S$ punctured at one point. Over the
punctured point the section will wind $a$ times
around the fiber. Such a section is a sub-complex of $G^*X$,
homeomorphic to a $2$-dimensional surface with boundary.
The integral of $\alpha$ over the boundary equals
$a$. Thus, according to the Stokes formula, the
integral of $d\alpha$ over the whole section also
equals $a$. But the integral of $d\alpha$ over the
section equals the integral of $\omega$
over $S$. \qed

\subsection{Back to the first Chern classes of $\L_i$}

The complex $A$ of stable ribbon graphs 
and the total space of the
polygonal bundle $\B$ are obviously
polytopal complexes. The projection $\B \rightarrow A$
is a morphism of complexes as in 
Proposition~\ref{Prop:Chern}.

It is straightforward to check that the $1$-form
$\alpha$ defined in Section~\ref{Ssec:conncurv} is
a $1$-form on the total space of $\B$ in the sense
of Definition~\ref{Def:forms}.

Thus it remains to check that $\alpha$ satisfies
the conditions of Proposition~\ref{Prop:Chern}.

\begin{proposition}{\rm \cite{Kontsevich}}

(i)~The integral of $\alpha$ over any fiber of $\B$ 
equals $-1$. 

(ii)~The $2$-from $d \alpha$ is the lifting 
of a $2$-form $\omega$ from the base $A$. 
\end{proposition}

\paragraph{Proof.}  \

(i)~As we go around the fiber,
the distinguished point goes around the polygon $B$
counterclockwise. The coefficients $l_i/p$ remain
constant, while each $\phi_i$ decreases from its
initial value to $0$ and then from $p$ back to its
initial value. Thus the integral over the fiber of
each $d (\phi_i/p)$ equals $-1$ and the sum of
the coefficients $l_i/p$ equals $1$.

(ii)~A simple calculation gives
$$
\omega = d\alpha = \sum_{1 \leq i < j \leq k-1} 
d \! \left( \frac{l_i}p \right)
\wedge
d \! \left( \frac{l_j}p \right)
$$
on each cell. This $2$-form depends only on the lengths
$l_i$ but not on the $\phi_i$s. Therefore it is a lifting 
of a $2$-from from the base $A$. \qed

\bigskip

Thus $\omega$ is a $2$-form on the polytopal complex $A$
that represents minus the first Chern class of the bundle $\B$.
To calculate the intersection numbers of these Chern classes
one can multiply the $2$-forms $\omega$ and integrate them
over the cells of highest dimension.

This finishes the proof of Theorem~\ref{Thm:int=int}.

\end{document}